\documentclass[12pt, reqno]{amsart}
\usepackage{amsmath, amsthm, amscd, amsfonts, amssymb, graphicx, color}
\usepackage[bookmarksnumbered, colorlinks, plainpages]{hyperref}

\hypersetup{colorlinks=true,linkcolor=red, anchorcolor=green, citecolor=cyan, urlcolor=red, filecolor=magenta, pdftoolbar=true}
\theoremstyle{definition}

\theoremstyle{remark}

\numberwithin{equation}{section}
\input{tcilatex}

\begin{document}
\title[Further results on strictly Lipschitz summing operators ]{Further
results on strictly Lipschitz summing operators }
\author{ Maatougui BELAALA and Khalil SAADI}
\date{}
\dedicatory{Laboratoire d'Analyse Fonctionnelle et G\'{e}om\'{e}trie des
Espaces, Universit\'{e} de M'sila, Alg\'{e}rie.\\
belaala.maato@gmail.com\\
kh\_saadi@yahoo.fr}

\begin{abstract}
We give some new characterizations of strictly Lipschitz $p$-summing
operators. These operators have been introduced in order to improve the
Lipschitz $p$-summing operators. Therefore, we adapt this definition for
constructing other classes of Lipschitz mappings which are called strictly
Lipschitz $p$-nuclear and strictly Lipschitz $\left( p,r,s\right) $-summing
operators. Some interesting properties and factorization results are
obtained for these new classes.
\end{abstract}

\maketitle

\setcounter{page}{1}


\let\thefootnote\relax\footnote{%
Copyright 2016 by the Tusi Mathematical Research Group.}

\section{Introduction and preliminaries}

In the last years, many concept of the theory of $p$-summing operators have
been developed in several ways, namely the multilinear and Lipschitz
settings. The fundamental purpose of nonlinear theory of Lipschitz mappings
is to attempt to borrow linear properties in order to make an analog in the
nonlinear case. Let $X$ be a pointed metric space and $E$ be a Banach space.
It is well known that each Lipschitz operator $T:X\rightarrow E$ can be
factorized through a Lipschitz map and a linear operator. Let $\mathcal{I}$
be an ideal linear, when we want to form an analogous class of Lipschitz
operators, it is natural to think how to preserve the connection between
their linearization operators, which appear in the factorization, and the
original ideal $\mathcal{I}$. If this property holds for a given Lipschitz
class, then it can be repesented by \cite{1} as follows 
\begin{equation}
\mathcal{I}\circ Lip_{0}\left( X;E\right) =\mathcal{I}\left( \mathcal{F}%
\left( X\right) ;E\right) .  \tag{1.1}
\end{equation}%
Then, the above representation, that we consider interesting, expresses good
relation between Lipschitz operators and their linearizations. To make the
relation (1.1) attainable, we have improved in \cite{16} the definition of
Lipschitz $p$-summing by introducing the strictly Lipschitz $p$-summing
operators whose original ideal is $\Pi _{p},$ the Banach space of $p$%
-summing operators, and admits a similar representation of (1.1). The goal
of this paper is to explore more properties of the class of strictly
Lipschitz $p$-summing by showing some characterizations of those operators
by means of fundamental inequalities and the domination theorem of Pietsch.
We give a strong version of this concept, called M-strictly Lipschitz $p$%
-summing, when we consider the class of Lipschitz operators defined on
metric spaces. Furthermore, we will introduce the class of strictly
Lipschitz $p$-nuclear by generalizing the linear definition introduced by
Cohen \cite{6}. Next, we will naturally define the class of strictly
Lipschitz $\left( p,r,s\right) $-summing. In analogy with the definition of
strictly Lipschitz $p$-nuclear, some interesting results are obtained for
the class of strictly Lipschitz $\left( p,r,s\right) $-summing.%
\vspace{0.5cm}%

The paper is organized as follows.

First, we recall some standard notations which will be used throughout this
paper. In section 2, we study some characterizations of strictly Lipschitz $p
$-summing operators which are defined from a metric space into a Banach
sapce. We will give the definition for a general case, called M-strictly
Lipschitz $p$-summing, where $X$ and $Y$ are metric spaces for which there
is a beautiful equivalence between $T$ and its corresponding linear
operators for the concept of $p$-summing. In Section 3, the definition of
strictly Lipschitz $p$-nuclear operators is given. This class has surprising
properties namely their connections with linearization operators. A
representation through a Lipschitz tensor product is given for this class.
We end this section by investigating certain relations with other classes.
Section 4 is devoted to study the class of strictly Lipschitz $\left(
p,r,s\right) $-summing operators. These operators involve Lipschitz $\left(
p,r,s\right) $-summing operators which are introduced by \cite{3}$.$ Many
considerations in the previous section are analogous to that in section 4.%
\vspace{0.5cm}%

Now, we recall briefly some basic notations and terminology which we need in
the sequel. Throughout this paper, the letters $E,F$ will denote Banach
spaces and $X,Y$ will denote metric spaces with a distinguished point
(pointed metric spaces) which we denote by $0$. Let $X$ be a pointed metric
space, we denote by $X^{\#}$ the Banach space of all Lipschitz functions $%
f:X\longrightarrow \mathbb{R}$ which vanish at $0$ under the Lipschitz norm
given by%
\begin{equation*}
Lip\left( f\right) =\sup \left\{ \frac{\left\vert f\left( x\right) -f\left(
y\right) \right\vert }{d\left( x,y\right) }:x,y\in X,x\neq y\right\} .
\end{equation*}%
We denote by $\mathcal{F}\left( X\right) ,$ the Lipschitz-free Banach space
over $X,$ the closed linear span of the linear forms $\delta _{\left(
x,y\right) }$ of $Lip_{0}\left( X\right) ^{\ast }$ such that 
\begin{equation*}
\delta _{\left( x,y\right) }\left( f\right) =f\left( x\right) -f\left(
y\right) ,\text{ for every }f\in Lip_{0}\left( X\right) ,
\end{equation*}%
i.e.,%
\begin{equation*}
\mathcal{F}\left( X\right) =\overline{span\left\{ \delta _{\left( x,y\right)
}:x,y\in X\right\} }^{Lip_{0}\left( X\right) ^{\ast }}
\end{equation*}%
We have $X^{\#}=\mathcal{F}\left( X\right) ^{\ast }$ holds isometrically via
the application%
\begin{equation*}
Q_{X}\left( f\right) \left( m\right) =m\left( f\right) ,\text{ for every }%
f\in X^{\#}\text{ and }m\in \mathcal{F}\left( X\right) .
\end{equation*}%
For the general theory of free Banach spaces, see \cite{9, 10, 14, 18}. Let $%
X$ be a pointed metric space and $E$ be a Banach space, we denote by $%
Lip_{0}\left( X;E\right) $ the Banach space of all Lipschitz functions
(Lipschitz operators) $T:X\rightarrow E$ such that $T\left( 0\right) =0$
with pointwise addition and Lipschitz norm. Note that for any $T\in
Lip_{0}\left( X;E\right) ,$ then there exists a unique linear map
(linearization of $T$) $\widehat{T}:\mathcal{F}\left( X\right)
\longrightarrow E$ such that $\widehat{T}\circ \delta _{X}=T$ and $%
\left\Vert \widehat{T}\right\Vert =Lip\left( T\right) ,$ i.e., the following
diagram commutes%
\begin{equation*}
\begin{array}{ccc}
X & \overset{T}{\longrightarrow } & Y \\ 
\delta _{X}\downarrow  & \nearrow \widehat{T} &  \\ 
\mathcal{F}\left( X\right)  &  & 
\end{array}%
\end{equation*}%
where $\delta _{X}$ is the canonical embedding so that $\left\langle \delta
_{X}\left( x\right) ,f\right\rangle =\delta _{\left( x,0\right) }\left(
f\right) =f\left( x\right) $ for $f\in X^{\#}.$ If $X$ is a Banach space and 
$T:X\rightarrow E$ is a linear operator, then the corresponding linear
operator $\widehat{T}$ is given by%
\begin{equation*}
\widehat{T}=T\circ \beta _{X},
\end{equation*}%
where $\beta _{X}:\mathcal{F}\left( X\right) \rightarrow X$ is linear
quotient map which verifies $\beta _{X}\circ \delta _{X}=id_{X}$ and $%
\left\Vert \beta _{X}\right\Vert \leq 1,$ see \cite[p 124]{10} for more
details about the operator $\beta _{X}$. Let $X,Y$ be two metric spaces. Let 
$T:X\rightarrow Y$ be a Lipschitz operator, then there is a unique linear
operator $\widetilde{T}$\ such that the following diagram commutes 
\begin{equation}
\begin{array}{ccc}
X & \underrightarrow{T} & Y \\ 
\downarrow \delta _{X} &  & \downarrow \delta _{Y} \\ 
F\left( X\right)  & \underrightarrow{\widetilde{T}} & F\left( Y\right) 
\end{array}
\tag{1.2}
\end{equation}%
i.e., $\delta _{Y}\circ T=\widetilde{T}\circ \delta _{X}$. The Lipschitz
adjoint map $T^{\#}:Y^{\#}\rightarrow X^{\#}$ of $T$ is defined as follows%
\begin{equation*}
T^{\#}\left( g\right) \left( x\right) =g\left( T\left( x\right) \right) ,%
\text{ for every }g\in Y^{\#}\text{ and }x\in X.
\end{equation*}%
The Lipschitz transpose operator $T_{/E^{\ast }}^{t}:E^{\ast }\rightarrow
X^{\#}$ is the restriction of $T^{\#}$ on $E^{\ast }$. We have, 
\begin{equation}
T^{\#}=Q_{X}^{-1}\circ \widetilde{T}^{\ast }\circ Q_{Y}\text{ and }%
T^{t}=Q_{X}^{-1}\circ \widehat{T}^{\ast }.  \tag{1.3}
\end{equation}%
Let $X$ be a metric space and $E$ be a Banach space, by $X\boxtimes E$ we
denote the Lipschitz tensor product of $X$ and $E$. This is the vector space
spanned by the linear functional $\delta _{\left( x,y\right) }\boxtimes e$
on $Lip_{0}\left( X;E^{\ast }\right) $ defined by%
\begin{equation*}
\delta _{\left( x,y\right) }\boxtimes e\left( f\right) =\left\langle f\left(
x\right) -f\left( y\right) ,e\right\rangle .
\end{equation*}%
See \cite{4} for more details about the properties of the space $X\boxtimes
E.$ Now, let $E$ be a Banach space, then $B_{E}$ denotes its closed unit
ball and $E^{\ast }$ its (topological) dual. Consider $1\leq p\leq \infty $
and $n\in \mathbb{N}^{\ast }$. We denote by $l_{p}^{n}\left( E\right) $ the
Banach space of all sequences $\left( e_{i}\right) _{i=1}^{n}$ in $E$ with
the norm

\begin{center}
$\left\Vert \left( e_{i}\right) _{i}\right\Vert _{l_{p}^{n}\left( E\right)
}=(\sum_{i=1}^{n}\left\Vert e_{i}\right\Vert ^{p})^{\frac{1}{p}}$,
\end{center}

\noindent and by $l_{p}^{n,w}\left( E\right) $ the Banach space of all
sequences $\left( e_{i}\right) _{i=1}^{n}$ in $E$ with the norm

\begin{center}
$\left\Vert \left( e_{i}\right) _{i}\right\Vert _{l_{p}^{n,w}\left( E\right)
}=\underset{e^{\ast }\in B_{E^{\ast }}}{\sup }(\sum_{i=1}^{n}\left\vert
\left\langle e_{i},e^{\ast }\right\rangle \right\vert ^{p})^{\frac{1}{p}}.$
\end{center}

\noindent If $E=\mathbb{K},$ we simply write $l_{p}^{n}$ and $l_{p}^{n,w}.%
\vspace{0.5cm}%
$

We recall some definitions which we need in the sequel. We refer to \cite{6,
7, 13}\ for more details about the following notions. Let\textit{\ }$1\leq
p\leq \infty $ and $p^{\ast }$ its conjugate, i.e., $\frac{1}{p}+\frac{1}{%
p^{\ast }}=1.$ Let $E,F$ be two Banach spaces and $R:E\rightarrow F$ be a
linear operator. Then,$%
\vspace{0.5cm}%
$

- The linear operator $R$\ is $p$-summing\ if there exists a constant $C>0$\
such that, for any $x_{1},...,x_{n}\in X,$ we have%
\begin{equation}
(\dsum\limits_{i=1}^{n}\left\Vert R\left( x_{i}\right) \right\Vert ^{p})^{%
\frac{1}{p}}\leq C\underset{x^{\ast }\in B_{X^{\ast }}}{\sup }%
(\tsum\limits_{i=1}^{n}\left\vert x^{\ast }\left( x_{i}\right) \right\vert
)^{\frac{1}{p}}.  \tag{1.4}
\end{equation}%
The class of $p$-summing linear operators from $E$ into $F$, which is
denoted by $\Pi _{p}(E;F),$ is a Banach space for the norm $\pi _{p}(R),$
i.e., the smallest constant $C$ such that the inequality (1.4) holds.%
\vspace{0.5cm}%

- The linear operator $R$\ is Cohen strongly $p$-summing\ if there exists a
constant $C>0$\ such that, for any $x_{1},...,x_{n}\in X,$\ and any $%
y_{1}^{\ast },...,y_{n}^{\ast }\in F^{\ast }$, we have%
\begin{equation}
\dsum\limits_{i=1}^{n}\left\vert \left\langle R\left( x_{i}\right)
,y_{i}^{\ast }\right\rangle \right\vert \leq C\left\Vert \left( x_{i}\right)
\right\Vert _{l_{p}^{n}\left( E\right) }\left\Vert \left( y_{i}^{\ast
}\right) \right\Vert _{l_{p^{\ast }}^{n,w}}.  \tag{1.5}
\end{equation}%
The class of Cohen strongly $p$-summing operators from $E$ into $F$, which
is denoted by $\mathcal{D}_{p}(E,F),$ is a Banach space for the norm $%
d_{p}(R),$ i.e., the smallest constant $C$ such that the inequality (1.5)
holds. If $1\leq p<\infty ,$ we have by \cite[Theorem 2.2.2]{6}%
\begin{equation}
R\in \Pi _{p}(E;F)\Leftrightarrow R^{\ast }\in \mathcal{D}_{p^{\ast
}}(F^{\ast };E^{\ast }).%
\vspace{0.5cm}
\tag{1.6}
\end{equation}

- The linear operator $R$\ is $p$-nuclear if there exists a constant $C>0$\
such that, for any $x_{1},...,x_{n}\in X,$\ and any $y_{1}^{\ast
},...,y_{n}^{\ast }\in F^{\ast }$, we have%
\begin{equation}
\dsum\limits_{i=1}^{n}\left\vert \left\langle R\left( x_{i}\right)
,y_{i}^{\ast }\right\rangle \right\vert \leq C\left\Vert \left( x_{i}\right)
\right\Vert _{l_{p}^{n,w}\left( E\right) }\left\Vert \left( y_{i}^{\ast
}\right) \right\Vert _{l_{p^{\ast }}^{n,w}}.  \tag{1.7}
\end{equation}%
The class of $p$-nuclear linear operators from $E$ into $F$, which is
denoted by $\mathcal{N}_{p}(E,F),$ is a Banach space for the norm $N_{p}(R),$
i.e., the smallest constant $C$ such that the inequality (1.7) holds. We
have by \cite{6, 12}: $R$ is $p$-nuclear if and only if, $R=R_{1}\circ R_{2}$
where $R_{1}$ is Cohen strongly $p$-summing and $R_{2}$ is $p$ summing.%
\vspace{0.5cm}%

- The linear operator $R$\ is $\left( p,r,s\right) $-summing if there exists
a constant $C>0$\ such that, for any $x_{1},...,x_{n}\in X,$\ and any $%
y_{1}^{\ast },...,y_{n}^{\ast }\in F^{\ast }$, we have%
\begin{equation}
\left\Vert \left( \left\langle R\left( x_{i}\right) ,y_{i}^{\ast
}\right\rangle \right) \right\Vert _{l_{p}^{n}}\leq C\left\Vert \left(
x_{i}\right) \right\Vert _{l_{r}^{n,w}\left( E\right) }\left\Vert \left(
y_{i}^{\ast }\right) \right\Vert _{l_{s}^{n,w}}.  \tag{1.8}
\end{equation}%
The class of $\left( p,r,s\right) $-summing linear operators from $E$ into $F
$, which is denoted by $\Pi _{p,r,s}(E,F),$ is a Banach space for the norm $%
\pi _{p,r,s}(R),$ i.e., the smallest constant $C$ such that the inequality
(1.8) holds. We have by \cite{13}: $R$ is $\left( p,r,s\right) $-summing if
and only if, $R=R_{1}\circ R_{2}$ where $R_{1}$ is Cohen strongly $s^{\ast }$%
-summing and $R_{2}$ is $r$-summing. In a particular case when $p=1$ and $%
\frac{1}{r}+\frac{1}{s}=1,$ we have%
\begin{equation*}
\mathcal{N}_{r}(E,F)=\Pi _{1,r,s}(E,F).
\end{equation*}

\section{\textsc{Characterization of strictly Lipschitz }$p$\textsc{-summing
operators}}

Let $X$ be a pointed metric space and $E$ be a Banach space. In order to
establish a relation between a Lipschitz operator $T:X\rightarrow E$ and its
linearization $\widehat{T}:\mathcal{F}\left( X\right) \rightarrow E$ for the
concept of $p$-summing, we have introduced in \cite{16} the notion of
strictly Lipschitz $p$-summing operators. Indeed, both operators are related
in the sense that: $T$ is strictly Lipschitz $p$-summing if and only if, $%
\widehat{T}$ is $p$-summing. In fact, that report was not true for the class
of Lipschitz $p$-summing operators which introduced by Farmer \cite{8}, see 
\cite[Remark 3.3]{15}. In this section, we give some characterizations of
the class of strictly Lipschitz $p$-summng operators and we adapt its
definition to the general case where the spaces are metric for which we
obtain a good relation between the Lipschitz operator $T$ and its
corresponding linear operators $\widetilde{T}$ and $T^{\#}$. Now, we start
by the definition of $d_{p}^{L},$ the corresponding Lipschitz cross-norm of
the tensor norm $d_{p}$, see \cite{16}\ fore more detail about this norm.
Recall the definition of the norms of Chevet-Saphar $d_{p}$ and $g_{p}$ \cite%
{5, 17} defined on Banach spaces, 
\begin{equation*}
d_{p}\left( u\right) =\inf \left\{ \left\Vert \left( x_{i}\right)
_{i}\right\Vert _{l_{p^{\ast }}^{n,w}\left( E\right) }\left\Vert \left(
y_{i}\right) _{i}\right\Vert _{l_{p}^{n}\left( F\right) }\right\} ,
\end{equation*}%
where the infimum is taking over all representations of $u$ of the form $%
u=\sum_{i=1}^{n}x_{i}\otimes y_{i}\in E\otimes F.$ The tensor norm $g_{p}$
is defined as follows%
\begin{equation*}
g_{p}(\sum_{i=1}^{n}x_{i}\otimes y_{i})=d_{p}^{t}(\sum_{i=1}^{n}x_{i}\otimes
y_{i})=d_{p}(\sum_{i=1}^{n}y_{i}\otimes x_{i})
\end{equation*}%
For every\textit{\ }$u=\sum_{k=1}^{l}\delta _{\left( x_{k},y_{k}\right)
}\boxtimes s_{k}\in X\boxtimes E,$\ we put%
\begin{equation}
A_{u}=\left\{ m=\sum_{i=1}^{n}m_{i}\otimes e_{i}\in \mathcal{F}\left(
X\right) \otimes E:m=\sum_{k=1}^{l}\delta _{\left( x_{k},y_{k}\right)
}\otimes s_{k}\right\} .  \tag{2.1}
\end{equation}%
Since the linearization $\widehat{T}$ can be seen as a linear form on $%
\mathcal{F}\left( X\right) \otimes E^{\ast },$ then for every $m\in A_{u}$
we have%
\begin{eqnarray}
\sum_{k=1}^{l}\left\langle T\left( x_{k}\right) -T\left( y_{k}\right)
,s_{k}^{\ast }\right\rangle  &=&\widehat{T}(\sum_{k=1}^{l}\delta _{\left(
x_{k},y_{k}\right) }\otimes s_{k}^{\ast })=\widehat{T}\left( m\right)  
\TCItag{2.2} \\
&=&\sum_{i=1}^{n_{1}}\left\langle \widehat{T}\left( m_{i}\right)
,e_{i}^{\ast }\right\rangle   \notag \\
&=&\sum_{i=1}^{n_{1}}\sum_{j=1}^{n_{2}}\left\langle \lambda _{i}^{j}(T\left(
x_{i}^{j}\right) -T\left( y_{i}^{j}\right) ),e_{i}^{\ast }\right\rangle . 
\notag
\end{eqnarray}%
where%
\begin{equation*}
m_{i}=\sum_{j=1}^{k_{i}}\lambda _{i}^{j}\delta _{\left(
x_{i}^{j},y_{i}^{j}\right) }=\sum_{j=1}^{n_{2}}\lambda _{i}^{j}\delta
_{\left( x_{i}^{j},y_{i}^{j}\right) },
\end{equation*}%
with $n_{2}=\max_{i=1}^{n_{1}}k_{i}$ and the terms between $k_{i}$ and $n_{2}
$ are zero. Now, Let $\alpha $ be a tensor norm defined on two Banach
spaces, by \cite[Theorem 3.1]{16}, there is a Lipschitz cross-norm $\alpha
^{L}$ which is defined on Lipschitz tensor product $X\boxtimes E$. Note that
if $u=\sum_{k=1}^{l}\delta _{\left( x_{k},y_{k}\right) }\boxtimes s_{k}\in
X\boxtimes E$ we have%
\begin{equation}
\alpha ^{L}(\sum_{k=1}^{l}\delta _{\left( x_{k},y_{k}\right) }\boxtimes
s_{k})=\alpha (\sum_{k=1}^{l}\delta _{\left( x_{k},y_{k}\right) }\otimes
s_{k}).  \tag{2.3}
\end{equation}%
where $\sum_{k=1}^{l}\delta _{\left( x_{k},y_{k}\right) }\otimes s_{k}\in 
\mathcal{F}\left( X\right) \otimes E$. So, we have%
\begin{equation*}
d_{p}^{L}\left( u\right) =\inf_{m\in A_{u}}\left\{ \left\Vert
m_{i}\right\Vert _{l_{p}^{n,w}\left( \mathcal{F}\left( X\right) \right)
}\left\Vert \left( e_{i}\right) _{i}\right\Vert _{l_{p^{\ast }}^{n}\left(
E\right) }\right\} .%
\vspace{0.5cm}%
\end{equation*}

\textbf{Definition 2.1 }\cite{16}. Let $1\leq p\leq \infty .$ A Lipschitz
operator $T:X\rightarrow E$ is said to be\textit{\ strictly Lipschitz }$p$%
\textit{-summing} if there exists a positive constant $C$ such that for
every $x_{k},y_{k}\in X$ and $s_{k}^{\ast }\in E^{\ast }$ $\left( 1\leq
k\leq l\right) $ we have 
\begin{equation}
\left\vert \sum_{k=1}^{l}\left\langle T\left( x_{k}\right) -T\left(
y_{k}\right) ,s_{k}^{\ast }\right\rangle \right\vert \leq Cd_{p}^{L}(u), 
\tag{2.4}
\end{equation}%
where $u=\sum_{k=1}^{l}\delta _{\left( x_{k},y_{k}\right) }\boxtimes
s_{k}^{\ast }.$ We denote by $\Pi _{p}^{SL}\left( X,E\right) $ the Banach
space of all strictly Lipschitz $p$-summing operators from $X$ into $E$
which its norm $\pi _{p}^{SL}\left( T\right) $ is the smallest constant $C$
verifying (2.4). If we consider linear operators defined on Banach spaces,
we have shown in \cite[Proposition 3.8]{16} that the three notions: $p$%
-summing, Lipschitz $p$-summing and strictly Lipschitz $p$-summing are\
coincide. The following characterization is the main result of this section.%
\vspace{0.5cm}%

\textbf{Theorem 2.2}. \textit{Let} $1\leq p\leq \infty .$ \textit{Let }$X$%
\textit{\ be a pointed metric space and }$E$\textit{\ be a Banach space. Let 
}$T:X\rightarrow E$\textit{\ be a Lipschitz operator. The following
properties are equivalent.}

\noindent \textit{1) }$T$\textit{\ is strictly Lipschitz }$p$\textit{%
-summing.}

\noindent \textit{2) }$\widehat{T}$\textit{\ is }$p$\textit{-summing.}

\noindent \textit{3) There exist a constant }$C>0$ \textit{and} \textit{a} 
\textit{Radon} \textit{probability }$\mu $\textit{\ on }$B_{X^{\#}}$ \textit{%
such that for all }$\left( x^{j}\right) _{j=1}^{n},\left( y^{j}\right)
_{j=1}^{n}\subset X$ and $\left( \lambda ^{j}\right) _{j=1}^{n}\subset 
\mathbb{K};\left( n\in \mathbb{N}^{\ast }\right) ,$ \textit{we have}%
\begin{equation}
\left\Vert \sum_{j=1}^{n}\lambda ^{j}\left( T\left( x^{j}\right) -T\left(
y^{j}\right) \right) \right\Vert \leq C(\dint\limits_{B_{X^{\#}}}\left\vert
\sum_{j=1}^{n}\lambda ^{j}\left( f\left( x^{j}\right) -f\left( y^{j}\right)
\right) \right\vert ^{p}d\mu \left( f\right) )^{\frac{1}{p}}.  \tag{2.5}
\end{equation}

\noindent \textit{4)} \textit{There is a constant }$C>0$\textit{\ such that
for every }$\left( x_{i}^{j}\right) _{i=1}^{n_{1}},\left( y_{i}^{j}\right)
_{i=1}^{n_{1}}$\textit{\ in }$X$ and $\left( \lambda _{i}^{j}\right)
_{i=1}^{n_{1}}\subset \mathbb{K};\left( 1\leq j\leq n_{2}\right) $ \textit{%
and} $n_{1},n_{2}\in \mathbb{N}^{\ast }$\textit{, we have}%
\begin{equation}
\sum_{i=1}^{n_{1}}(\left\Vert \sum_{j=1}^{n_{2}}\lambda _{i}^{j}\left(
T(x_{i}^{j})-T(y_{i}^{j})\right) \right\Vert ^{p})^{\frac{1}{p}}\leq
C\sup_{f\in X^{\#}}(\sum_{i=1}^{n_{1}}\left\Vert \sum_{j=1}^{n_{2}}\lambda
_{i}^{j}\left( f(x_{i}^{j})-f(y_{i}^{j})\right) \right\Vert ^{p})^{\frac{1}{p%
}}.%
\vspace{0.5cm}
\tag{2.6}
\end{equation}

\textit{Proof. }$\left( 1\right) \Longrightarrow \left( 2\right) :$ Theorem
3.5 in \cite{16}$.$

\noindent $\left( 2\right) \Longrightarrow \left( 3\right) :$ We apply
Pietsch Domination Theorem for $p$-summing linear operators \cite[Theorem
2.12]{7}$,$ then there is a Radon probability $\mu $\ on\textit{\ }$%
B_{X^{\#}}$ such that for any $m\in \mathcal{F}\left( X\right) $ we have 
\begin{equation*}
\left\Vert \widehat{T}\left( m\right) \right\Vert \leq
C(\dint\limits_{B_{X^{\#}}}\left\vert f\left( m\right) \right\vert ^{p}d\mu
\left( f\right) )^{\frac{1}{p}}.
\end{equation*}%
Now, let $\left( x^{j}\right) _{j=1}^{n},\left( y^{j}\right)
_{j=1}^{n}\subset X$ and $\left( \lambda ^{j}\right) _{j=1}^{n}\subset 
\mathbb{K},$ we put%
\begin{equation*}
m=\sum_{j=1}^{n}\lambda ^{j}\delta _{\left( x^{j},y^{j}\right) }\in \mathcal{%
F}\left( X\right) ,
\end{equation*}%
Then%
\begin{equation*}
\left\Vert \widehat{T}(\sum_{j=1}^{n}\lambda ^{j}\delta _{\left(
x^{j},y^{j}\right) })\right\Vert \leq C(\dint\limits_{B_{X^{\#}}}\left\vert
f(\sum_{j=1}^{n}\lambda ^{j}\delta _{\left( x^{j},y^{j}\right) })\right\vert
^{p}d\mu \left( f\right) )^{\frac{1}{p}}
\end{equation*}%
thus%
\begin{equation*}
\left\Vert \sum_{j=1}^{n}\lambda ^{j}\left( T\left( x^{j}\right) -T\left(
y^{j}\right) \right) \right\Vert \leq C(\dint\limits_{B_{X^{\#}}}\left\vert
\sum_{j=1}^{n}\lambda ^{j}\left( f\left( x^{j}\right) -f\left( y^{j}\right)
\right) \right\vert ^{p}d\mu \left( f\right) )^{\frac{1}{p}}.
\end{equation*}

$\left( 3\right) \Longrightarrow \left( 4\right) :$ Let $\left(
x_{i}^{j}\right) _{i=1}^{n_{1}},\left( y_{i}^{j}\right) _{i=1}^{n_{1}}$%
\textit{\ }in\textit{\ }$X$ and $\left( \lambda _{i}^{j}\right)
_{i=1}^{n_{1}}\subset \mathbb{K}\left( 1\leq j\leq n_{2}\right) $, by (2.5)
we have for every $1\leq i\leq n_{1}$%
\begin{equation*}
\left\Vert \sum_{j=1}^{n_{2}}\lambda _{i}^{j}\left(
T(x_{i}^{j})-T(y_{i}^{j})\right) \right\Vert \leq
C(\dint\limits_{B_{X^{\#}}}\left\vert \sum_{j=1}^{n_{2}}\lambda
_{i}^{j}\left( f(x_{i}^{j})-f(y_{i}^{j})\right) \right\vert ^{p}d\mu \left(
f\right) )^{\frac{1}{p}}
\end{equation*}%
Therefore,

\QTP{Body Math}
$\sum_{i=1}^{n_{1}}\left\Vert \sum_{j=1}^{n_{2}}\lambda _{i}^{j}(T\left(
x_{i}^{j}\right) -T\left( y_{i}^{j}\right) )\right\Vert ^{p}$

\QTP{Body Math}
$\leq C\dint\limits_{B_{X^{\#}}}\sum_{i=1}^{n_{1}}\left\vert
\sum_{j=1}^{n_{2}}\lambda _{i}^{j}(f\left( x_{i}^{j}\right) -f\left(
y_{i}^{j}\right) )\right\vert ^{p}d\mu \left( f\right) $

\QTP{Body Math}
$\leq \sup_{f\in B_{X^{\#}}}\sum_{i=1}^{n_{1}}\left\vert
\sum_{j=1}^{n_{2}}\lambda _{i}^{j}(f\left( x_{i}^{j}\right) -f\left(
y_{i}^{j}\right) )\right\vert ^{p}.$

\noindent Finally, we have 
\begin{equation*}
(\sum_{i=1}^{n_{1}}\left\Vert \sum_{j=1}^{n_{2}}\lambda _{i}^{j}\left(
T(x_{i}^{j})-T(y_{i}^{j})\right) \right\Vert ^{p})^{\frac{1}{p}}\leq
C\sup_{f\in X^{\#}}(\sum_{i=1}^{n_{1}}\left\vert \sum_{j=1}^{n_{2}}\lambda
_{i}^{j}\left( f(x_{i}^{j})-f(y_{i}^{j})\right) \right\vert ^{p})^{\frac{1}{p%
}}
\end{equation*}

$\left( 4\right) \Longrightarrow \left( 1\right) :$ Let $u=\sum_{k=1}^{l}%
\delta _{\left( x_{k},y_{k}\right) }\boxtimes s_{k}^{\ast }\in X\boxtimes
E^{\ast }$. Let $A_{u}$ the set as defined in (2.1). Let $%
m=\sum_{i=1}^{n_{1}}m_{i}\otimes e_{i}^{\ast }\in
A_{u}(m_{i}=\sum_{j=1}^{n_{2}}\lambda _{i}^{j}\delta _{\left(
x_{i}^{j},y_{i}^{j}\right) }),$ by (2.2)

\QTP{Body Math}
$\left\vert \sum_{k=1}^{l}\left\langle T\left( x_{k}\right) -T\left(
y_{k}\right) ,s_{k}^{\ast }\right\rangle \right\vert $

\QTP{Body Math}
$=\left\vert \sum_{i=1}^{n_{1}}\sum_{j=1}^{n_{2}}\left\langle \lambda
_{i}^{j}\left( T\left( x_{i}^{j}\right) -T\left( y_{i}^{j}\right) \right)
,e_{i}^{\ast }\right\rangle \right\vert $ by H\"{o}lder

\QTP{Body Math}
$\leq (\sum_{i=1}^{n_{1}}\left\Vert \sum_{j=1}^{n_{2}}\lambda _{i}^{j}\left(
T\left( x_{i}^{j}\right) -T\left( y_{i}^{j}\right) \right) \right\Vert
^{p})^{\frac{1}{p}}(\sum_{i=1}^{n_{1}}\left\Vert e_{i}^{\ast }\right\Vert
^{p^{\ast }})^{\frac{1}{p^{\ast }}}$

\QTP{Body Math}
$\leq C\sup_{f\in X^{\#}}(\sum_{i=1}^{n_{1}}\left\vert
\sum_{j=1}^{n_{2}}\lambda _{i}^{j}\left( f\left( x_{i}^{j}\right) -f\left(
y_{i}^{j}\right) \right) \right\vert ^{p})^{\frac{1}{p}}(\sum_{i=1}^{n_{1}}%
\left\Vert e_{i}^{\ast }\right\Vert ^{p^{\ast }})^{\frac{1}{p^{\ast }}}$

\QTP{Body Math}
$\leq C\left\Vert \left( m_{i}\right) \right\Vert _{l_{p}^{n_{1},w}\left( 
\mathcal{F}\left( X\right) \right) }\left\Vert \left( e_{i}^{\ast }\right)
\right\Vert _{l_{p^{\ast }}^{n_{1}}\left( E\right) }.$

\noindent By taking the infimum over all representations of $m\in A_{u},$ we
obtain%
\begin{equation*}
\left\vert \sum_{k=1}^{l}\left\langle T\left( x_{k}\right) -T\left(
y_{k}\right) ,s_{k}^{\ast }\right\rangle \right\vert \leq Cd_{p}^{L}(u).
\end{equation*}%
Then, $T$ is strictly Lipschitz $p$-summing.$\quad \blacksquare 
\vspace{0.5cm}%
$

If we put $n_{2}=1$ in the formula (2.6), we obtain exactly the definition
of Lipschitz $p$-summing operators. Now, let $X$ and $Y$ be two metric
spaces, Farmer \cite{8} has introduced the definition of Lipschitz $p$%
-summing for the general case of metric spaces. Inspired by the definition
of strictly Lipschitz $p$-summing, we give the following general setting.$%
\vspace{0.5cm}%
$

\textbf{Definition 2.3}. Let $1\leq p\leq \infty .$ Let $X$ and $Y$ be two
metric spaces. A Lipschitz operator $T:X\rightarrow Y$ is said to be\textit{%
\ M-strictly Lipschitz }$p$\textit{-summing} if there exists a positive
constant $C$ such that for every $x_{k},y_{k}\in X$ and $g_{k}\in Y^{\#}$ $%
\left( 1\leq k\leq l\right) $ we have 
\begin{equation*}
\left\vert \sum_{k=1}^{l}g_{k}\left( T\left( x_{k}\right) \right)
-g_{k}\left( T\left( y_{k}\right) \right) \right\vert \leq Cd_{p}^{L}(u),
\end{equation*}%
where $u=\sum_{k=1}^{l}\delta _{\left( x_{k},y_{k}\right) }\boxtimes
g_{k}\in X\boxtimes Y^{\#}.$ Note that if $Y$ is a Banach space, for every $%
y^{\ast }\in Y^{\ast }$ we have $Lip\left( y^{\ast }\right) =\left\Vert
y^{\ast }\right\Vert ,$ then%
\begin{equation*}
\left\Vert \left( y_{i}^{\ast }\right) \right\Vert _{l_{p^{\ast }}^{n}\left(
Y^{\#}\right) }=\left\Vert \left( y_{i}^{\ast }\right) \right\Vert
_{l_{p^{\ast }}^{n}\left( Y^{\ast }\right) }.
\end{equation*}%
Therefore, the above definition leads to that of strictly Lipschitz $p$%
-summing. As an interesting characterization of M-strictly Lipschitz $p$%
-summing operators, we have the following result.$%
\vspace{0.5cm}%
$

\textbf{Proposition 2.4}. \textit{Let }$1\leq p\leq \infty .$\textit{\ Let }$%
X$ \textit{and} $Y$\textit{\ be two metric spaces. The following properties
are equivalent.}

\noindent \textit{1) }$T:X\rightarrow Y$\textit{\ is M-strictly Lipschitz }$%
p $\textit{-summing.}

\noindent \textit{2) The linearization operator }$\widetilde{T}:\mathcal{F}%
\left( X\right) \rightarrow \mathcal{F}\left( Y\right) $\textit{\ is }$p$%
\textit{-summing.}

\noindent \textit{3) The Lipschitz adjoint }$T^{\#}:Y^{\#}\rightarrow X^{\#}$%
\textit{\ is Cohen strongly }$p$\textit{-summing}$.%
\vspace{0.5cm}%
$

For the proof, we need the following Lemma.$%
\vspace{0.5cm}%
$

\textbf{Lemma 2.5}. \textit{Let }$X$\textit{\ be a pointed metric space and }%
$E,F$\textit{\ be two Banach spaces. Suppose that }$E$\textit{\ and }$F$ 
\textit{are isometrically isomorphic via the application }$Q.$\textit{\ Let }%
$u=\sum_{k=1}^{l}\delta _{\left( x_{k},y_{k}\right) }\boxtimes s_{k}\in
X\boxtimes E$\textit{, then}%
\begin{equation*}
d_{p}^{L}\left( u\right) =d_{p}^{L}\left( w\right) ,
\end{equation*}%
\textit{where} $w=\sum_{k=1}^{l}\delta _{\left( x_{k},y_{k}\right)
}\boxtimes Q\left( s_{k}\right) \in X\boxtimes F.%
\vspace{0.5cm}%
$

\textit{Proof}. The identification $\mathcal{F}\left( X\right) \otimes E=%
\mathcal{F}\left( X\right) \otimes F$ holds via the transformation%
\begin{equation*}
\sum_{i=1}^{n}m_{i}\otimes e_{i}\mapsto \sum_{i=1}^{n}m_{i}\otimes Q\left(
e_{i}\right) .
\end{equation*}%
Let $m=\sum_{i=1}^{n}m_{i}\otimes e_{i}\in A_{u},$ then $\sum_{i=1}^{n}m_{i}%
\otimes Q\left( e_{i}\right) \in A_{w}.$ We have%
\begin{eqnarray*}
d_{p}^{L}\left( w\right) &\leq &\left\Vert m_{i}\right\Vert
_{l_{p}^{n,w}\left( \mathcal{F}\left( X\right) \right) }\left\Vert \left(
Q\left( e_{i}\right) \right) _{i}\right\Vert _{l_{p^{\ast }}^{n}\left(
F\right) } \\
&\leq &\left\Vert m_{i}\right\Vert _{l_{p}^{n,w}\left( \mathcal{F}\left(
X\right) \right) }\left\Vert \left( e_{i}\right) _{i}\right\Vert
_{l_{p^{\ast }}^{n}\left( E\right) },
\end{eqnarray*}%
by taking the infimum on $A_{u}$, we find $d_{p}^{L}\left( w\right) \leq
d_{p}^{L}\left( u\right) .$ With the same argument, the inverse inequality
is immediate.$\quad \blacksquare 
\vspace{0.5cm}%
$

\textit{Proof of Proposition 2.4}. $\left( 1\right) \Rightarrow \left(
2\right) :$ Let $\left( x_{k}\right) _{k=1}^{l}\subset X$ and $\left(
f_{k}^{\ast }\right) _{k=1}^{l}\subset \mathcal{F}\left( Y\right) ^{\ast }.$
We will show that $\delta _{Y}\circ T$ is strictly Lipschitz $p$-summing.
So, there is $\left( g_{k}\right) _{k=1}^{l}\subset Y^{\#}$ such that $%
f_{k}^{\ast }=Q_{Y}\left( g_{k}\right) $ for every $1\leq k\leq l.$ We have%
\begin{eqnarray*}
\left\vert \sum_{k=1}^{l}\left\langle \delta _{Y}\circ T\left( x_{k}\right)
-\delta _{Y}\circ T\left( y_{k}\right) ,f_{k}^{\ast }\right\rangle
\right\vert  &=&\left\vert \sum_{k=1}^{l}\left\langle \delta _{\left(
T\left( x_{k}\right) ,T\left( y_{k}\right) \right) },Q_{Y}\left(
g_{k}\right) \right\rangle \right\vert  \\
&=&\left\vert \sum_{k=1}^{l}g_{k}\left( T\left( x_{k}\right) \right)
-g_{k}\left( T\left( y_{k}\right) \right) \right\vert  \\
&\leq &Cd_{p}^{L}(u)=Cd_{p}^{L}(w),
\end{eqnarray*}%
where $u=\sum_{k=1}^{l}\delta _{\left( x_{k},y_{k}\right) }\boxtimes
g_{k}\in X\boxtimes Y^{\#}$ and $w=\sum_{k=1}^{l}\delta _{\left(
x_{k},y_{k}\right) }\boxtimes f_{k}^{\ast }\in X\boxtimes \mathcal{F}\left(
Y\right) ^{\ast }.$ Then, $\delta _{Y}\circ T$ is strictly Lipschitz $p$%
-summing. Now, by the factorization (1.2) we have%
\begin{equation*}
\delta _{Y}\circ T=\widetilde{T}\circ \delta _{X},
\end{equation*}%
this shows that the linearization of $\delta _{Y}\circ T$ is $\widetilde{T}.$
Then, by Theorem 2.2, $\widetilde{T}$ is $p$-summing.

$\left( 2\right) \Rightarrow \left( 3\right) :$ We have by (1.3)%
\begin{equation*}
T^{\#}=Q_{X}^{-1}\circ \widetilde{T}^{\ast }\circ Q_{Y},
\end{equation*}%
hence, by (1.6) and the ideal property, $T^{\#}$ is Cohen strongly $p^{\ast }
$-summing.

$\left( 3\right) \Rightarrow \left( 1\right) :$ Suppose that $T^{\#}$ is
Cohen strongly $p^{\ast }$-summing. Let $u=\sum_{k=1}^{l}\delta _{\left(
x_{k},y_{k}\right) }\boxtimes g_{k}\in X\boxtimes Y^{\#}$ then

\QTP{Body Math}
$\left\vert \sum_{k=1}^{l}g_{k}\left( T\left( x_{k}\right) \right)
-g_{k}\left( T\left( y_{k}\right) \right) \right\vert $

\QTP{Body Math}
$=\left\vert \sum_{k=1}^{l}\left\langle T^{\#}\left( g_{k}\right) ,\delta
_{\left( x_{k},y_{k}\right) }\right\rangle \right\vert $

\QTP{Body Math}
$\leq Cg_{p}(\sum_{k=1}^{l}g_{k}\otimes \delta _{\left( x_{k},y_{k}\right)
}) $ (we know that $g_{p}=d_{p}^{t}$)

\QTP{Body Math}
$\leq Cd_{p}(\sum_{k=1}^{l}\delta _{\left( x_{k},y_{k}\right) }\otimes
g_{k})($by (2.3)$)$

\QTP{Body Math}
$=Cd_{p}^{L}\left( u\right) .$

\noindent Therefore, $T$ is M-strictly Lipschitz $p$-summing.$\quad
\blacksquare 
\vspace{0.5cm}%
$

The definition of Lipschitz dual of a given operator ideal $\mathcal{I}$\ is
given in \cite{1}\ as follows%
\begin{equation*}
\mathcal{I}^{Lip_{0}-dual}\left( X,E\right) =\left\{ T\in Lip_{0}\left(
X,E\right) :T^{t}\in \mathcal{I}\left( E^{\ast },X^{\#}\right) \right\} .%
\vspace{0.5cm}%
\end{equation*}

\textbf{Corollary 2.6}. \textit{Let }$\mathcal{D}_{p}$\textit{\ be the
linear ideal of Cohen strongly }$p$-\textit{summing operators. Then, we have}%
\begin{equation*}
\mathcal{D}_{p}^{Lip_{0}-dual}\left( X,E\right) =\Pi _{p^{\ast }}^{SL}\left(
X,E\right) .%
\vspace{0.5cm}%
\end{equation*}

\textit{Proof}. Let $T\in \Pi _{p^{\ast }}^{SL}\left( X,E\right) $, by
Theorem 2.2 its linearization $\widehat{T}$ is $p^{\ast }$-summing. By (1.3)
and the ideal property, $T^{t}$ is Cohen strongly $p$-summing, then $T\in 
\mathcal{D}_{p}^{Lip_{0}-dual}\left( X,E\right) .$ Conversely, let $T\in 
\mathcal{D}_{p}^{Lip_{0}-dual}\left( X,E\right) $, then $T^{t}:E^{\ast
}\rightarrow X^{\#}$ is Cohen strongly $p$-summing. We have%
\begin{equation*}
\widehat{T}^{\ast }=Q_{X}\circ T^{t},
\end{equation*}%
then $\widehat{T}$ is $p^{\ast }$-summing by\ (1.6) and the result follows
by Theorem 2.2.$\quad \blacksquare $

\section{\textsc{Strictly Lipschitz }$p$\textsc{-nuclear operators}}

In this section, we adopt the same procedure of the previous section for
defining the strictly Lipschitz $p$-nuclear operators. Then, we obtain a
nice class that has many interesting properties with good relations with
other classes of Lipschitz mappings where the linearization operators play
the key of all obtained results. Among the results of this section, we
obtain an analog of Pietsch domination theorem and then the same
factorization result to Cohen-Kwapie\'{n} which is given for $p$-nuclear
linear operators. Let $p\in \left[ 1,\infty \right] \ $and $E,F$ be two
Banach spaces. The definition of tensor norm $w_{p}$ on $E\otimes F$ is
given in \cite{6} by%
\begin{equation*}
w_{p}\left( u\right) =\inf \left\{ \left\Vert \left( x_{i}\right)
\right\Vert _{l_{p}^{n,w}\left( E\right) }\left\Vert \left( y_{i}\right)
\right\Vert _{l_{p^{\ast }}^{n,w}\left( F\right) }\right\} ,
\end{equation*}%
where the infimum is taken over all representations of $u$ of the form $%
u=\sum_{i=1}^{n}x_{i}\otimes y_{i}\in E\otimes F.$ Then by \cite[Lemma 2.5.1]%
{6}, we have the following identification%
\begin{equation}
\mathcal{N}_{p}\left( E,F\right) =\left( E\widehat{\otimes }_{w_{p}}F^{\ast
}\right) ^{\ast }.  \tag{3.1}
\end{equation}

Next, we will give an analog approach for the class of strictly Lipschitz $p$%
-nuclear operators. Let $X$ be a pointed metric space and $E$ be a Banach
space. Let $u=\sum_{k=1}^{l}\delta _{\left( x_{k},y_{k}\right) }\boxtimes
s_{k}^{\ast }\in X\boxtimes E^{\ast }$ and $A_{u}$ be the set as in (2.1).
We consider%
\begin{equation*}
w_{p}^{SL}\left( u\right) =\inf_{m\in A_{u}}\left\{ \left\Vert \left(
m_{i}\right) _{i}\right\Vert _{l_{p}^{n,w}\left( \mathcal{F}\left( X\right)
\right) }\left\Vert \left( e_{i}^{\ast }\right) \right\Vert _{l_{p^{\ast
}}^{n,w}\left( E^{\ast }\right) }\right\} .
\end{equation*}

The following Proposition shows that $w_{p}^{SL}$\ is nothing else than the
corresponding Lipschitz cross-norm to the tensor norm $w_{p}.$%
\vspace{0.5cm}%

\textbf{Proposition 3.1}. \textit{Let }$X$\textit{\ be a pointed metric
space and }$E$\textit{\ be a Banach space. For every }$u\in X\boxtimes E$%
\textit{\ we have}%
\begin{equation*}
w_{p}^{SL}\left( u\right) =w_{p}^{L}\left( u\right) ,
\end{equation*}%
\textit{where }$w_{p}^{L}$\textit{\ is the Lipschitz cross-norm
corresponding to the tensor norm }$w_{p}.$%
\vspace{0.5cm}%

\textit{Proof}. Let $u=\sum_{k=1}^{l}\delta _{\left( x_{k},y_{k}\right)
}\boxtimes s_{k}\in X\boxtimes E,$ we have

\QTP{Body Math}
$w_{p}^{SL}\left( u\right) =\inf_{m\in A_{u}}\left\{ \left\Vert \left(
m_{i}\right) \right\Vert _{l_{p}^{n,w}\left( \mathcal{F}\left( X\right)
\right) }\left\Vert \left( e_{i}\right) \right\Vert _{l_{p^{\ast
}}^{n,w}\left( E\right) }\right\} $

\QTP{Body Math}
$=w_{p}(\sum_{k=1}^{l}\delta _{\left( x_{k},y_{k}\right) }\otimes
s_{k})\left( \text{by (2.3)}\right) $

\QTP{Body Math}
$=w_{p}^{L}\left( u\right) .\quad \blacksquare 
\vspace{0.5cm}%
$

It is not difficult to prove the following results.$%
\vspace{0.5cm}%
$

\textbf{Proposition 3.2}. \textit{Let }$X$\textit{\ be a pointed metric
space and }$E$\textit{\ be a Banach space. }

\noindent $\left( 1\right) $ \textit{For every }$u\in X\boxtimes E$\textit{\
we have }$w_{p}^{SL}\left( u\right) \leq d_{p}^{L}\left( u\right) .$

\noindent $\left( 2\right) $ \textit{If }$p=1$\textit{, we have }$%
w_{1}^{SL}\left( u\right) =d_{1}^{L}\left( u\right) .%
\vspace{0.5cm}%
$

Now, we give the following definition of strictly Lipschitz $p$-nuclear
operators for which we use the Lipschitz cross-norms.$%
\vspace{0.5cm}%
$

\textbf{Definition 3.3}. Let $1\leq p\leq \infty .$ The Lipschitz operator $%
T:X\rightarrow E$\ is strictly Lipschitz $p$-nuclear if\ for every $%
x_{k},y_{k}\in X$\ and $s_{k}^{\ast }\in E^{\ast }$\ $\left( 1\leq k\leq
l\right) $\ we have\textit{\ }%
\begin{equation}
\left\vert \sum_{k=1}^{l}\left\langle T\left( x_{k}\right) -T\left(
y_{k}\right) ,s_{k}^{\ast }\right\rangle \right\vert \leq Cw_{p}^{SL}(u), 
\tag{3.2}
\end{equation}%
where $u=\sum_{k=1}^{l}\delta _{\left( x_{k},y_{k}\right) }\boxtimes
s_{k}^{\ast }.$ The class of all strictly Lipschitz $p$-nuclear operators
from $X$ into $E$ is denoted by $\mathcal{N}_{p}^{SL}\left( X,E\right) $,
which is a Banach space with the norm $N_{p}^{SL}(T)$ which is the smallest
constant $C$ such that the inequality (3.2) holds.%
\vspace{0.5cm}%

As an immediate consequence of the Proposition 3.2, the following assertions
can be checked easily.%
\vspace{0.5cm}%

\textbf{Proposition 3.4}. \textit{Let }$X$\textit{\ be a pointed metric
space and }$E$\textit{\ be a Banach space.}

\noindent $\left( 1\right) $ \textit{For every }$1<p\leq \infty $\textit{\
we have}%
\begin{equation*}
\mathcal{N}_{p}^{SL}\left( X,E\right) \subset \Pi _{p}^{SL}\left( X,E\right)
.
\end{equation*}

\noindent $\left( 2\right) $ \textit{If }$p=1$\textit{, then}%
\begin{equation*}
\mathcal{N}_{1}^{SL}\left( X,E\right) =\Pi _{1}^{SL}\left( X,E\right) .%
\vspace{0.5cm}%
\end{equation*}

In the following result, we connect between a Lipschitz operator and its
linearization for the concept of $p$-nuclear. So, this connection helps in
considering the space of these operators as a Lipschitz ideal generated by
the composition method from the linear ideal $\mathcal{N}_{p}$.%
\vspace{0.5cm}%

\textbf{Theorem 3.5}. \textit{Let} $1\leq p\leq \infty .$ \textit{Let }$X$%
\textit{\ be a pointed metric space and }$E$\textit{\ be a Banach space. Let 
}$T:X\rightarrow E$\textit{\ be a Lipschitz operator. The following
properties are equivalent.}

\noindent \textit{1) }$T$\textit{\ is strictly Lipschitz }$p$\textit{%
-nuclear.}

\noindent \textit{2)} \textit{The linearization operator }$\widehat{T}$%
\textit{\ is }$p$\textit{-nuclear.}

\noindent \textit{On the other hand}%
\begin{equation*}
\mathcal{N}_{p}^{SL}\left( X,E\right) =\mathcal{N}_{p}\left( \mathcal{F}%
\left( X\right) ,E\right) \text{\textit{\ holds isometrically}.}%
\vspace{0.5cm}%
\end{equation*}

\textit{Proof}. Suppose that $\widehat{T}$ is $p$-nuclear. Let $x_{k},y_{k}$
in $X$ and $e_{k}^{\ast }\in E^{\ast }\left( 1\leq k\leq l\right) ,$ we put $%
u=\sum_{k=1}^{l}\delta _{\left( x_{k},y_{k}\right) }\boxtimes s_{k}^{\ast }.$
Let $A_{u}$ the set as defined in (1.2). For every $m\in A_{u},$ by (2.2) we
have%
\begin{equation*}
\sum_{k=1}^{l}\left\langle T\left( x_{k}\right) -T\left( y_{k}\right)
,s_{k}^{\ast }\right\rangle =\sum_{i=1}^{n_{1}}\left\langle \widehat{T}%
\left( m_{i}\right) ,e_{i}^{\ast }\right\rangle .
\end{equation*}%
Now, let us prove that $T$ is strictly Lipschitz $p$-nuclear. Let $%
m=\sum_{i=1}^{n_{1}}m_{i}\otimes e_{i}^{\ast }\in A_{u},$ by (1.7) we have%
\begin{eqnarray*}
\left\vert \sum_{k=1}^{l}\left\langle T\left( x_{k}\right) -T\left(
y_{k}\right) ,s_{k}^{\ast }\right\rangle \right\vert &=&\left\vert
\sum_{i=1}^{n_{1}}\left\langle \widehat{T}\left( m_{i}\right) ,e_{i}^{\ast
}\right\rangle \right\vert \\
&\leq &N_{p}\left( \widehat{T}\right) \left\Vert \left( m_{i}\right)
\right\Vert _{l_{p}^{n_{1},w}\left( \mathcal{F}\left( X\right) \right)
}\left\Vert \left( e_{i}^{\ast }\right) \right\Vert _{l_{p^{\ast
}}^{n_{1},w}\left( E\right) }.
\end{eqnarray*}%
By taking the infimum over all representations of $m\in A_{u},$ we obtain%
\begin{eqnarray*}
\left\vert \sum_{k=1}^{l}\left\langle T\left( x_{k}\right) -T\left(
y_{k}\right) ,s_{k}^{\ast }\right\rangle \right\vert &\leq &N_{p}\left( 
\widehat{T}\right) w_{p}(\sum_{k=1}^{l}\delta _{\left( x_{k},y_{k}\right)
}\otimes s_{k}^{\ast }) \\
&\leq &N_{p}\left( \widehat{T}\right) w_{p}^{SL}\left( u\right) .
\end{eqnarray*}%
Then, $T$ is strictly Lipschitz $p$-nuclear and%
\begin{equation*}
N_{p}^{SL}\left( T\right) \leq N_{p}\left( \widehat{T}\right) .
\end{equation*}%
Conversely, let $m_{i}\in \mathcal{F}\left( X\right) $ and $e_{i}^{\ast }\in
E^{\ast }\left( 1\leq i\leq n_{1}\right) $. Then

$\sum_{i=1}^{n_{1}}\left\vert \left\langle \widehat{T}\left( m_{i}\right)
,e_{i}^{\ast }\right\rangle \right\vert $

$=\sup_{\left( \xi _{i}\right) _{i}\in B_{l_{\infty }^{n_{1}}}}\left\vert
\sum_{i=1}^{n_{1}}\xi _{i}\left\langle \widehat{T}\left( m_{i}\right)
,e_{i}^{\ast }\right\rangle \right\vert $

$=\sup_{\left( \xi _{i}\right) _{i}\in B_{l_{\infty }^{n_{1}}}}\left\vert
\sum_{i=1}^{n_{1}}\xi _{i}\left\langle \sum_{j=1}^{n_{2}}\lambda
_{i}^{j}\left( T\left( x_{i}^{j}\right) -T\left( y_{i}^{j}\right) \right)
,e_{i}^{\ast }\right\rangle \right\vert $

$=\sup_{\left( \xi _{i}\right) _{i}\in B_{l_{\infty }^{n_{1}}}}\left\vert
\left\langle \sum_{i=1}^{n_{1}}\sum_{j=1}^{n_{2}}T\left( x_{i}^{j}\right)
-T\left( y_{i}^{j}\right) ,\lambda _{i}^{j}\xi _{i}e_{i}^{\ast
}\right\rangle \right\vert $

$\leq \sup_{\left( \xi _{i}\right) _{i}\in B_{l_{\infty
}^{n_{1}}}}N_{p}^{SL}\left( T\right) w_{p}^{SL}\left( u\right) ,$

\noindent where 
\begin{eqnarray*}
u &=&\sum_{i=1}^{n_{1}}\sum_{j=1}^{n_{2}}\delta _{\left(
x_{i}^{j},y_{i}^{j}\right) }\boxtimes \lambda _{i}^{j}\xi _{i}e_{i}^{\ast }
\\
&=&\sum_{i=1}^{n_{1}}\xi _{i}\sum_{j=1}^{n_{2}}\lambda _{i}^{j}\delta
_{\left( x_{i}^{j},y_{i}^{j}\right) }\boxtimes e_{i}^{\ast } \\
&=&\sum_{i=1}^{n_{1}}\xi _{i}m_{i}\boxtimes e_{i}^{\ast }
\end{eqnarray*}%
Therefore,%
\begin{eqnarray*}
\sum_{i=1}^{n_{1}}\left\vert \left\langle \widehat{T}\left( m_{i}\right)
,e_{i}^{\ast }\right\rangle \right\vert &\leq &\sup_{\left( \xi _{i}\right)
_{i}\in B_{l_{\infty }^{n_{1}}}}N_{p}^{SL}\left( T\right) \left\Vert \left(
\xi _{i}m_{i}\right) \right\Vert _{l_{p^{\ast }}^{n_{1},w}\left( \mathcal{F}%
\left( X\right) \right) }\left\Vert \left( e_{i}^{\ast }\right) \right\Vert
_{l_{p^{\ast }}^{n_{1},w}\left( E^{\ast }\right) } \\
&\leq &N_{p}^{SL}\left( T\right) \left\Vert \left( m_{i}\right) \right\Vert
_{l_{p^{\ast }}^{n_{1},w}\left( \mathcal{F}\left( X\right) \right)
}\left\Vert \left( e_{i}^{\ast }\right) \right\Vert _{l_{p^{\ast
}}^{n_{1},w}\left( E^{\ast }\right) }.
\end{eqnarray*}%
Then, by (1.7) $\widehat{T}$ is $p$-nuclear and 
\begin{equation*}
N_{p}\left( \widehat{T}\right) \leq N_{p}^{SL}\left( T\right) .\quad
\blacksquare 
\vspace{0.5cm}%
\end{equation*}

The Pietsch domination theorem is one of the interesting characterizations
which is verified by the class of strictly Lipschitz $p$-nuclear. It can be
proved by means of the linearization operators. Consequently, there is an
equivalent definition in which we use only fundamental inequalities.$%
\vspace{0.5cm}%
$

\textbf{Theorem 3.6}. \textit{Let }$X$\textit{\ be a pointed metric space
and }$E$\textit{\ be a Banach space. Let }$T:X\rightarrow E$\textit{\ be a
Lipschitz operator. The following properties are equivalent.}

\noindent \textit{1) }$T$\textit{\ is strictly Lipschitz }$p$\textit{%
-nuclear.}

\noindent \textit{2) There exist a constant }$C>0$, \textit{a} \textit{Radon}
\textit{probability }$\mu $\textit{\ on }$B_{X^{\#}}$ \textit{and} $\eta \in
B_{E^{\ast \ast }}$ \textit{such that for every }$\left( x^{j}\right)
_{j=1}^{n},\left( y^{j}\right) _{j=1}^{n}\subset X,$ $\left( \lambda
^{j}\right) _{j=1}^{n}\subset \mathbb{K}$ \textit{and} $e^{\ast }\in E^{\ast
},$ \textit{we have}%
\begin{equation}
\left\vert \left\langle \sum_{j=1}^{n}\lambda ^{j}\left( T\left(
x^{j}\right) -T\left( y^{j}\right) \right) ,e^{\ast }\right\rangle
\right\vert \leq C(\dint\limits_{B_{X^{\#}}}\left\vert \sum_{j=1}^{n}\lambda
^{j}\left( f\left( x^{j}\right) -f\left( y^{j}\right) \right) \right\vert
^{p}d\mu \left( f\right) )^{\frac{1}{p}}\left\Vert e^{\ast }\right\Vert
_{L_{p^{\ast }}\left( \eta \right) }.  \tag{3.3}
\end{equation}

\noindent \textit{3) For any} $\left( x_{i}^{j}\right) _{i=1}^{n_{1}},\left(
y_{i}^{j}\right) _{i=1}^{n_{1}}\subset X_{j},$ $\left( \lambda
_{i}^{j}\right) _{i=1}^{n_{1}}\subset \mathbb{K},$\ ($j=1,...,n_{2}$)\ 
\textit{and any }$\left( e_{i}^{\ast }\right) _{i=1}^{n_{1}}\subset E^{\ast }
$\textit{, we have}%
\begin{equation}
\sum_{i=1}^{n_{1}}\left\vert \left\langle \sum_{j=1}^{n_{2}}\lambda
_{i}^{j}(T(x_{i}^{j})-T(y_{i}^{j})),e_{i}^{\ast }\right\rangle \right\vert
\leq C\sup_{f\in B_{X^{\#}}}(\sum_{i=1}^{n_{1}}\left\vert
\sum_{j=1}^{n_{2}}\lambda _{i}^{j}(f(x_{i}^{j})-f(y_{i}^{j}))\right\vert
^{p})^{\frac{1}{p}}\left\Vert (e_{i}^{\ast })\right\Vert _{l_{p^{\ast
}}^{n_{1}w}}.  \tag{3.4}
\end{equation}%
\textit{In this case, we have}%
\begin{eqnarray*}
N_{p}^{SL}\left( T\right)  &=&\inf \left\{ C:\text{ \textit{verifying (3.3)}}%
\right\}  \\
&=&\inf \left\{ C:\text{ \textit{verifying (3.4)}}\right\} .%
\vspace{0.5cm}%
\end{eqnarray*}

\textit{Proof}. $\left( 1\right) \Rightarrow \left( 2\right) :$ Let $\left(
x^{j}\right) _{j=1}^{n},\left( y^{j}\right) _{j=1}^{n}\subset X,\left(
\lambda ^{j}\right) _{j=1}^{n}\subset \mathbb{K}$ and $e^{\ast }\in E^{\ast
}.$ We put%
\begin{equation*}
m=\sum_{j=1}^{n}\lambda ^{j}\delta _{\left( x^{j},y^{j}\right) },
\end{equation*}%
By Theroem 3.5, the linearization operator $\widehat{T}$ is $p$-nuclear.
Then, by \cite[Proposition 2]{12}, there exist a Radon probability $\mu $\
on $B_{X^{\#}}$ and $\eta \in B_{E^{\ast \ast }}$ such that

\QTP{Body Math}
$\left\vert \left\langle \sum_{j=1}^{n}\lambda ^{j}\left( T\left(
x^{j}\right) -T\left( y^{j}\right) \right) ,e^{\ast }\right\rangle
\right\vert =\left\vert \left\langle \widehat{T}\left( m\right) ,e^{\ast
}\right\rangle \right\vert $

\QTP{Body Math}
$\leq (\dint\limits_{B_{X^{\#}}}\left\vert f\left( m\right) \right\vert
^{p}d\mu \left( f\right) )^{\frac{1}{p}}\left\Vert e^{\ast }\right\Vert
_{L_{p^{\ast }}\left( \eta \right) }$

\QTP{Body Math}
$\leq C(\dint\limits_{B_{X^{\#}}}\left\vert \sum_{j=1}^{n}\lambda ^{j}\left(
f\left( x^{j}\right) -f\left( y^{j}\right) \right) \right\vert ^{p}d\mu
\left( f\right) )^{\frac{1}{p}}\left\Vert e^{\ast }\right\Vert _{L_{p^{\ast
}}\left( \eta \right) }.$

$\left( 2\right) \Rightarrow \left( 3\right) :$ Let $\left( x_{i}^{j}\right)
_{i=1}^{n_{1}},\left( y_{i}^{j}\right) _{i=1}^{n_{1}}\subset X,\left(
\lambda _{i}^{j}\right) _{i=1}^{n_{1}}\subset \mathbb{K}\left( 1\leq j\leq
n_{2}\right) $ and $e_{1}^{\ast },...,e_{n_{1}}^{\ast }\in E^{\ast }.$ By
(3.3) we have for every $1\leq i\leq n_{1}$%
\begin{equation*}
\left\vert \left\langle \sum_{j=1}^{n_{2}}\lambda _{i}^{j}\left(
T(x_{i}^{j})-T(y_{i}^{j})\right) ,e_{i}^{\ast }\right\rangle \right\vert
\leq C(\dint\limits_{B_{X^{\#}}}\left\vert \sum_{j=1}^{n_{2}}\lambda
_{i}^{j}\left( f(x_{i}^{j})-f(y_{i}^{j})\right) \right\vert ^{p}d\mu \left(
f\right) )^{\frac{1}{p}}\left\Vert e_{i}^{\ast }\right\Vert _{L_{p^{\ast
}}\left( \eta \right) }.
\end{equation*}%
So, we have

\QTP{Body Math}
$\sum_{i=1}^{n_{1}}\left\vert \left\langle \sum_{j=1}^{n_{2}}\lambda
_{i}^{j}\left( T\left( x_{i}^{j}\right) -T\left( y_{i}^{j}\right) \right)
,e_{i}^{\ast }\right\rangle \right\vert $

\QTP{Body Math}
$\leq C\sum_{i=1}^{n_{1}}(\dint\limits_{B_{X^{\#}}}\left\vert
\sum_{j=1}^{n_{2}}\lambda _{i}^{j}\left( f(x_{i}^{j})-f(y_{i}^{j})\right)
\right\vert ^{p}d\mu \left( f\right) )^{\frac{1}{p}}\left\Vert e_{i}^{\ast
}\right\Vert _{L_{p^{\ast }}\left( \eta \right) }$ by H\"{o}lder

\QTP{Body Math}
$\leq C(\dint\limits_{B_{X^{\#}}}\sum_{i=1}^{n_{1}}\left\vert
\sum_{j=1}^{n_{2}}\lambda _{i}^{j}\left( f(x_{i}^{j})-f(y_{i}^{j})\right)
\right\vert ^{p}d\mu \left( f\right) )^{\frac{1}{p}}(\dint\limits_{B_{E^{%
\ast \ast }}}\sum_{i=1}^{n_{1}}\left\vert e^{\ast \ast }\left( e_{i}^{\ast
}\right) \right\vert ^{p^{\ast }}d\eta \left( e^{\ast \ast }\right) )^{\frac{%
1}{p^{\ast }}}$

\QTP{Body Math}
$\leq C\sup_{f\in B_{X^{\#}}}(\sum_{i=1}^{n_{1}}\left\vert
\sum_{j=1}^{n_{2}}\lambda _{i}^{j}\left( f\left( x_{i}^{j}\right) -f\left(
y_{i}^{j}\right) \right) \right\vert ^{p})^{\frac{1}{p}}\left\Vert \left(
e_{i}^{\ast }\right) \right\Vert _{l_{p^{\ast }}^{n_{1},w}\left( E^{\ast
}\right) }.$

$\left( 3\right) \Rightarrow \left( 1\right) :$ Let $x_{k},y_{k}$ in $X$ and 
$e_{k}^{\ast }\in E^{\ast }\left( 1\leq k\leq l\right) ,$ we put $%
u=\sum_{k=1}^{l}\delta _{\left( x_{k},y_{k}\right) }\boxtimes s_{k}^{\ast }.$
Let $m=\sum_{i=1}^{n_{1}}m_{i}\otimes e_{i}^{\ast }\in A_{u}$. By (2.6)%
\begin{eqnarray*}
\left\vert \sum_{k=1}^{l}\left\langle T\left( x_{k}\right) -T\left(
y_{k}\right) ,s_{k}^{\ast }\right\rangle \right\vert &=&\left\vert
\sum_{i=1}^{n_{1}}\left\langle \sum_{j=1}^{n_{2}}\lambda _{i}^{j}\left(
T\left( x_{i}^{j}\right) -T\left( y_{i}^{j}\right) \right) ,e_{i}^{\ast
}\right\rangle \right\vert \\
&\leq &C\sup_{f\in B_{X^{\#}}}(\sum_{i=1}^{n_{1}}\left\vert
\sum_{j=1}^{n_{2}}\lambda _{i}^{j}\left( f\left( x_{i}^{j}\right) -f\left(
y_{i}^{j}\right) \right) \right\vert ^{p})^{\frac{1}{p}}\left\Vert \left(
e_{i}^{\ast }\right) \right\Vert _{l_{p^{\ast }}^{n_{1},w}\left( E^{\ast
}\right) } \\
&\leq &C\left\Vert \left( m_{i}\right) \right\Vert _{l_{p^{\ast
}}^{n_{1},w}\left( \mathcal{F}\left( X\right) \right) }\left\Vert \left(
e_{i}^{\ast }\right) \right\Vert _{l_{p^{\ast }}^{n_{1},w}\left( E^{\ast
}\right) }.
\end{eqnarray*}%
By taking the infimum over all representations of $m\in A_{u},$ we obtain%
\begin{equation*}
\left\vert \sum_{k=1}^{l}\left\langle T\left( x_{k}\right) -T\left(
y_{k}\right) ,s_{k}^{\ast }\right\rangle \right\vert \leq Cw_{p}^{SL}\left(
u\right) .
\end{equation*}%
Then, $T$ is strictly Lipschitz $p$-nuclear.$\quad \blacksquare 
\vspace{0.5cm}%
$

In the following proposition, we give new examples of strictly Lipschitz $p$%
-nuclear operators.%
\vspace{0.5cm}%

\textbf{Proposition 3.7.} \textit{Let }$1\leq p\leq \infty .$\textit{\ Let }$%
X$\textit{\ be a pointed metric space, }$E$\textit{\ and }$F$\textit{\ be
Banach spaces. Let }$v\in \mathcal{B}\left( E,F\right) $\textit{\ and }$T\in
Lip_{0}\left( X,E\right) .$

\noindent $\left( 1\right) $\textit{\ If }$T$\textit{\ is strictly Lipschitz 
}$p$\textit{-nuclear, then }$v\circ T$ \textit{strictly Lipschitz }$p$%
\textit{-nuclear. We have }%
\begin{equation*}
N_{p}^{SL}(v\circ T)\leq \left\Vert v\right\Vert N_{p}^{SL}(T).
\end{equation*}%
$\left( 2\right) $ \textit{If }$L\in \Pi _{p}^{SL}\left( X,E\right) $\textit{%
\ and }$v\in \mathcal{D}_{p}\left( E,F\right) $\textit{, then }$v\circ L$ 
\textit{is strictly Lipschitz }$p$\textit{-nuclear and }%
\begin{equation*}
N_{p}^{SL}(v\circ L)\leq d_{p}\left( v\right) \pi _{p}^{SL}\left( L\right) .%
\vspace{0.5cm}%
\end{equation*}

\textbf{Theorem 3.8}. \textit{Let }$X$\textit{\ be a pointed metric space
and }$E$\textit{\ be a Banach space. Let }$T:X\rightarrow E$\textit{\ be a
Lipschitz operator. The following properties are equivalent.}

\noindent \textit{1) }$T$\textit{\ is strictly Lipschitz }$p$\textit{%
-nuclear.}

\noindent \textit{2) There exist a Banach space }$G,$\textit{\ a strictly
Lipschitz }$p$\textit{-summing operator }$R:X\rightarrow G$\textit{\ and a
Cohen strongly }$p$\textit{-summing linear operator }$S:G\rightarrow E$%
\textit{\ such that }$T=S\circ R.%
\vspace{0.5cm}%
$

\textit{Proof.} The second implication is immediate by Proposition 3.7. For
the first, let $T$\textit{\ }be a\textit{\ }strictly Lipschitz $p$-nuclear
operator. Since $\widehat{T}$ is $p$-nuclear, it factors as follows

\begin{equation*}
\begin{array}{ccccc}
\widehat{T}: & \mathcal{F}\left( X\right) & \rightarrow &  & E \\ 
& \searrow &  & \nearrow &  \\ 
& L & G & S & 
\end{array}%
\end{equation*}%
where $L$ is $p$-summing and $S$ is Cohen strongly $p$-summing linear
operators. Consequently, $T$ admits the following factorization%
\begin{equation*}
\begin{array}{cccc}
T: & X & \rightarrow & E \\ 
& \delta _{X}\downarrow &  & \uparrow S \\ 
& \mathcal{F}\left( X\right) & \underrightarrow{L} & G%
\end{array}%
\end{equation*}%
On the other hand, $T=S\circ R$ where $R=L\circ \delta _{X}.$ So, as $%
\widehat{R}=L$ which is $p$-summing, $R$ is strictly Lipschitz $p$-summing
by Theorem 2.2.$\quad \blacksquare 
\vspace{0.5cm}%
$

According to the last results and by applying \cite[Corollary 3.2]{16}, we
can identify the space of strictly Lipschitz $p$-nuclear operators with the
dual of the space $X\boxtimes E$ endowing with the norm $w_{p}^{SL}.$ This
identification views as a Lipschitz analog of (3.1).$%
\vspace{0.5cm}%
$

\textbf{Corollary 3.9}. \textit{Let }$X$\textit{\ be a pointed metric space
and }$E$\textit{\ be a Banach space. Then}%
\begin{equation*}
\mathcal{N}_{p}^{SL}\left( X,E\right) =(X\widehat{\boxtimes }%
_{w_{p}^{SL}}E^{\ast })^{\ast }=(\mathcal{F}\left( X\right) \widehat{\otimes 
}_{w_{p}}E^{\ast })^{\ast }.
\end{equation*}

Now, we survey some relations between certain classes of Lipschitz operators
basing on linear results and by using the linearization operators. First, in
the linear case we know that a $p$-nuclear operator is compact whenever its
domain or range is reflexive. In \cite{11}, there is an equivalence between
a Lipschitz map and its linearization, so we have the following fact.%
\vspace{0.5cm}%

\textbf{Corollary 3.10.} \textit{Let }$1<p\leq \infty .$\textit{\ Let }$X$%
\textit{\ be a pointed metric space and }$E$\textit{\ be a reflexive Banach
space. Then, every strictly Lipschitz }$p$\textit{-nuclear operator is
compact.}%
\vspace{0.5cm}%

\textbf{Corollary 3.11.} \textit{Let }$1<p\leq \infty .$\textit{\ Let }$X$%
\textit{\ be a pointed metric space and }$E$ \textit{be Banach space. Then }

\noindent \textit{(1) }$\mathcal{N}_{p}^{SL}\left( X,E\right) \subset 
\mathcal{D}_{p}^{L}\left( X,E\right) .$

\noindent \textit{(2) }$\mathcal{N}_{p}^{SL}\left( X,E\right) \subset \Pi
_{p}^{SL}\left( X,E\right) \subset \Pi _{p}^{L}\left( X,E\right) .$%
\vspace{0.5cm}%

Next, we have the same linear inclusion and a good relation with Lipschitz
G-integral operators which introduced in \cite{2}. By \cite[Proposition 2.4]%
{2}, \cite[Theorem 2.5.2]{6} and Theorem 3.5, we have the following
inclusion.%
\vspace{0.5cm}%

\textbf{Corollary 3.12.} \textit{Let }$1<p\leq \infty .$\textit{\ Let }$X$%
\textit{\ be a pointed metric space and }$E$ \textit{be Banach space. Then}%
\begin{equation*}
Lip_{0GI}\left( X,E\right) \subset N_{p}^{SL}\left( X,E\right) .%
\vspace{0.5cm}%
\end{equation*}

If $E$ is an $\mathcal{L}_{p^{\ast }}$-space, we obtain by \cite[Theorem
3.3.3]{6} the following coincidence.%
\vspace{0.5cm}%

\textbf{Corollary 3.13.} \textit{Let }$1<p\leq \infty .$\textit{\ Let }$X$%
\textit{\ be a pointed metric space and }$E$ \textit{is an} $\mathcal{L}%
_{p^{\ast }}$\textit{-space. Then, the Lipschitz G-integral operators and
strictly Lipschitz }$p$\textit{-nuclear operators are coincide.}%
\vspace{0.5cm}%

As a particular case, every Hilbert space $H$ is an $\mathcal{L}_{2}$-space,
we have%
\begin{equation*}
N_{2}^{SL}\left( X,H\right) =Lip_{0GI}\left( X,H\right) .%
\vspace{0.5cm}%
\end{equation*}

\textbf{Corollary 3.14}. \textit{Let }$X,E$\textit{\ be two Banach spaces
and }$\mathcal{N}_{p}$\textit{\ the linear ideal of }$p$-\textit{nuclear
operators. Then we have,}%
\begin{equation*}
\mathcal{N}_{p}^{Lip_{0}-dual}\left( X,E\right) =\mathcal{N}_{p^{\ast
}}^{SL}\left( X,E\right) .
\end{equation*}

\textit{Proof}. Let $T\in \mathcal{N}_{p^{\ast }}^{SL}\left( X,E\right) $,
then by Theorem 3.5, its linearization $\widehat{T}$ is $p^{\ast }$-nuclear.
By the ideal property and (1.3), $T^{t}$ is $p$-nuclear, then%
\begin{equation*}
\mathcal{N}_{p^{\ast }}^{SL}\left( X,E\right) \subset \mathcal{N}%
_{p}^{Lip_{0}-dual}\left( X,E\right) .
\end{equation*}%
Let $T\in \mathcal{N}_{p}^{Lip_{0}-dual}\left( X,E\right) $, then $%
T^{t}:E^{\ast }\rightarrow X^{\#}$ is $p$-nuclear, then%
\begin{equation*}
T^{t}=v_{1}\circ v_{2}
\end{equation*}%
where $v_{1}$ is Cohen strongly $p$-summing and $v_{2}$ is $p$-summing.
Therefore, $Q_{X}^{-1}\circ \widehat{T}^{\ast }=v_{1}\circ v_{2},$ and then 
\begin{equation*}
\widehat{T}^{\ast }=Q_{X}\circ v_{1}\circ v_{2},
\end{equation*}%
which is $p$-nuclear by the ideal property. So, $\widehat{T}$ is $p^{\ast }$%
-nuclear by \cite[Theorem 2.2.4]{6}. Finally, $T$ is strictly Lipschitz $%
p^{\ast }$-nuclear.$\quad \blacksquare 
\vspace{0.5cm}%
$

\section{\textsc{Strictly Lipschitz }$\left( p,r,s\right) $\textsc{-summing
operators}}

In the same circle of ideas, we study the strong version of Lipschitz $%
\left( p,r,s\right) $-summing linear operators. The linear class has been
stated\ by Lapreste in \cite{13} and generalized to Lipschitz case by Ch\'{a}%
vez-Dom\'{\i}nguez \cite{3}. The results of this section are analogous of
the setting of strictly Lipschitz $p$-nuclear operators. Now, we recall the
following definition as stated in \cite{3}.$%
\vspace{0.5cm}%
$

\textbf{Definition 4.1}. Let $X$ be a pointed metric space and $E$ be a
Banach space. Let $T:X\rightarrow E$ be a Lipschitz map. $T$ is Lipschitz $%
\left( p,r,s\right) $-summing if there is a constant $C>0$\ such that for
any $n\in \mathbb{N}^{\ast }$, $\left( x_{i}\right) _{i},\left( y_{i}\right)
_{i}$ in $X;\left( e_{i}^{\ast }\right) _{i}$ in $Y^{\ast }$ and $\left(
\lambda _{i}\right) _{i},\left( k_{i}\right) _{i}$ in $\mathbb{R}_{+}^{\ast }
$ $\left( 1\leq i\leq n\right) $, we have%
\begin{equation}
\left\Vert \left( \lambda _{i}\left\langle T\left( x_{i}\right) -T\left(
y_{i}\right) ,e_{i}^{\ast }\right\rangle \right) _{i}\right\Vert
_{l_{p}^{n}}\leq Cw_{r}^{Lip}\left( \left( \lambda
_{i}k_{i}^{-1},x_{i},y_{i}\right) _{i}\right) \left\Vert \left(
k_{i}e_{i}^{\ast }\right) _{i}\right\Vert _{l_{s}^{n,w}\left( E^{\ast
}\right) },  \tag{4.1}
\end{equation}%
where $w_{r}^{Lip}\left( \left( \lambda _{i}k_{i}^{-1},x_{i},e_{i}\right)
_{i=1}^{n}\right) $ is the weak Lipschitz $p$-norm defined by%
\begin{eqnarray*}
w_{r}^{Lip}\left( \left( \lambda _{i},x_{i},y_{i}\right) _{i=1}^{n}\right) 
&=&\sup_{f\in B_{X^{\#}}}(\sum_{i=1}^{n}\left\vert \lambda _{i}\left(
f\left( x_{i}\right) -f\left( y_{i}\right) \right) \right\vert ^{r})^{\frac{1%
}{r}} \\
&=&\left\Vert \left( \lambda _{i}\delta _{\left( x_{i},y_{i}\right) }\right)
\right\Vert _{l_{r}^{n,w}\left( \mathcal{F}\left( X\right) \right) }.
\end{eqnarray*}%
We denote by $\Pi _{p,r,s}^{L}\left( X,E\right) $ the Banach space of all
Lipschitz $\left( p,r,s\right) $-summing operators with the norm $\pi
_{p,r,s}^{L}(T)$ which is the smallest constant $C$ such that the inequality
(4.1) holds.%
\vspace{0.5cm}%

\textbf{Definition 4.2.} Let $0<p,r,s<\infty $ such that $\frac{1}{p}\leq 
\frac{1}{r}+\frac{1}{s}.$ Let $X$ be a pointed metric space and $E$ be a
Banach space. The Lipschitz operator $T:X\rightarrow E$ is strictly
Lipschitz $\left( p,r,s\right) $-summing if there is a constant $C>0$\ such
that for any $n\in \mathbb{N}^{\ast }$, $\left( x_{i}^{j}\right)
_{i=1}^{n_{1}},\left( y_{i}^{j}\right) _{i=1}^{n_{1}}\subset X,$\ $\left(
\lambda _{i}^{j}\right) _{i=1}^{n_{1}}\subset \mathbb{K}$($j=1,...,n_{2}$)\
and any $e_{1}^{\ast },...,e_{n_{1}}^{\ast }\in E^{\ast }$, we have%
\begin{equation}
(\dsum\limits_{i=1}^{n_{1}}\left\vert \left\langle
\dsum\limits_{j=1}^{n_{2}}\lambda
_{i}^{j}(T(x_{i}^{j})-T(y_{i}^{j})),e_{i}^{\ast }\right\rangle \right\vert
^{p})^{\frac{1}{p}}\leq C\sup_{f\in
B_{X^{\#}}}(\dsum\limits_{i=1}^{n_{1}}\left\vert
\dsum\limits_{j=1}^{n_{2}}\lambda
_{i}^{j}(f(x_{i}^{j})-f(y_{i}^{j}))\right\vert ^{r})^{\frac{1}{r}}\left\Vert
(e_{i}^{\ast })\right\Vert _{l_{s}^{n_{1}w}}.  \tag{4.2}
\end{equation}%
The class of all strictly Lipschitz $\left( p,r,s\right) $-summing operators
from $X$ into $E$ is denoted by $\Pi _{p,r,s}^{SL}\left( X,E\right) $, which
is a Banach space with the norm $\pi _{p,r,s}^{SL}(T)$ which is the smallest
constant $C$ such that the inequality (4.2) holds.%
\vspace{0.5cm}%

\textbf{Remark 4.3. }$\left( 1\right) $ If we put $n_{2}=1$, we obtain 
\begin{equation*}
(\dsum\limits_{i=1}^{n_{1}}\left\vert \left\langle \lambda _{i}\left(
T(x_{i})-T(y_{i})\right) ,e_{i}^{\ast }\right\rangle \right\vert ^{p})^{%
\frac{1}{p}}\leq C\sup_{f\in
B_{X^{\#}}}(\dsum\limits_{i=1}^{n_{1}}\left\vert k_{i}^{-1}\lambda
_{i}\left( f\left( x_{i}\right) -f\left( y_{i}\right) \right) \right\vert
^{r})^{\frac{1}{r}}\left\Vert (k_{i}e_{i}^{\ast })\right\Vert
_{l_{s}^{n_{1}w}}
\end{equation*}%
i.e., 
\begin{equation*}
\Pi _{p,r,s}^{SL}\left( X,E\right) \subset \Pi _{p,r,s}^{L}\left( X,E\right)
.
\end{equation*}

\noindent $\left( 2\right) $ If $p=1,$ $r$ and $s$ verify $\frac{1}{r}+\frac{%
1}{s}=1$, the definition coincides with the definition of strictly Lipschitz 
$r$-nuclear.%
\vspace{0.5cm}%

\textbf{Theorem 4.4}. \textit{Suppose that }$\frac{1}{p}=\frac{1}{r}+\frac{1%
}{s}.$\textit{\ Let }$X$\textit{\ be a pointed metric space and }$E$\textit{%
\ be a Banach space. Let }$T:X\rightarrow E$\textit{\ be a Lipschitz
operator. The following properties are equivalent.}

\noindent \textit{1) }$T$\textit{\ is strictly Lipschitz }$\left(
p,r,s\right) $\textit{-summing.}

\noindent \textit{2)} $\widehat{T}$ \textit{is }$\left( p,r,s\right) $%
\textit{-summing.}

\noindent \textit{In this case we have} $\pi _{p,r,s}\left( \widehat{T}%
\right) =\pi _{p,r,s}^{SL}\left( T\right) .$%
\vspace{0.5cm}%

\textit{Proof}. $\left( 1\right) \Rightarrow \left( 2\right) :$ Let $T$ be a
strictly Lipschitz $\left( p,r,s\right) $-summing operator. Let $m_{i}\in 
\mathcal{F}\left( X\right) $ of the form $m_{i}=\sum_{j=1}^{n_{2}}\lambda
_{i}^{j}\delta _{\left( x_{i}^{j},y_{i}^{j}\right) }$ and $e_{1}^{\ast
},...,e_{n_{1}}^{\ast }\in E^{\ast }.$ Then%
\begin{eqnarray*}
(\dsum\limits_{i=1}^{n_{1}}\left\vert \left\langle \widehat{T}\left(
m_{i}\right) ,e_{i}^{\ast }\right\rangle \right\vert ^{p})^{\frac{1}{p}}
&=&(\dsum\limits_{i=1}^{n_{1}}\left\vert \left\langle
\dsum\limits_{j=1}^{n_{2}}\lambda _{i}^{j}\left(
T(x_{i}^{j})-T(y_{i}^{j})\right) ,e_{i}^{\ast }\right\rangle \right\vert
^{p})^{\frac{1}{p}} \\
&\leq &\pi _{p,r,s}^{SL}\left( T\right) \sup_{f\in
B_{X^{\#}}}(\dsum\limits_{i=1}^{n_{1}}\left\vert f\left( m_{i}\right)
\right\vert ^{r})^{\frac{1}{r}}\left\Vert (e_{i}^{\ast })\right\Vert
_{l_{s}^{n_{1}w}} \\
&\leq &\pi _{p,r,s}^{SL}\left( T\right) \left\Vert (m_{i})\right\Vert
_{l_{r}^{n_{1},w}\left( \mathcal{F}\left( X\right) \right) }\left\Vert
(e_{i}^{\ast })\right\Vert _{l_{s}^{n_{1}w}}.
\end{eqnarray*}%
Then, $\widehat{T}$ is\textit{\ }$\left( p,r,s\right) $-summing and we have%
\begin{equation*}
\pi _{p,r,s}\left( \widehat{T}\right) \leq \pi _{p,r,s}^{SL}\left( T\right) .
\end{equation*}%
$\left( 2\right) \Rightarrow \left( 1\right) :$ Let $\left( x_{i}^{j}\right)
_{i=1}^{n_{1}},\left( y_{i}^{j}\right) _{i=1}^{n_{1}}\in X,$\ $\left(
\lambda _{i}^{j}\right) _{i=1}^{n_{1}}\subset \mathbb{K}$($j=1,...,n_{2}$)\
and any $e_{1}^{\ast },...,e_{n_{1}}^{\ast }\in E^{\ast }$

\QTP{Body Math}
$(\dsum\limits_{i=1}^{n_{1}}\left\vert \left\langle
\dsum\limits_{j=1}^{n_{2}}\lambda _{i}^{j}\left(
T(x_{i}^{j})-T(y_{i}^{j})\right) ,e_{i}^{\ast }\right\rangle \right\vert
^{p})^{\frac{1}{p}}$

\QTP{Body Math}
$=(\dsum\limits_{i=1}^{n_{1}}\left\vert \left\langle \widehat{T}\left(
m_{i}\right) ,e_{i}^{\ast }\right\rangle \right\vert ^{p})^{\frac{1}{p}}$

\QTP{Body Math}
$\leq \pi _{p,r,s}\left( \widehat{T}\right) \left\Vert (m_{i})\right\Vert
_{l_{r}^{n_{1},w}\left( \mathcal{F}\left( X\right) \right) }\left\Vert
(e_{i}^{\ast })\right\Vert _{l_{s}^{n_{1}w}}$

\QTP{Body Math}
$\leq \pi _{p,r,s}\left( \widehat{T}\right) \sup_{f\in
B_{X^{\#}}}(\dsum\limits_{i=1}^{n_{1}}\left\vert
\dsum\limits_{j=1}^{n_{2}}\lambda _{i}^{j}\left(
f(x_{i}^{j})-f(y_{i}^{j})\right) \right\vert ^{r})^{\frac{1}{r}}\left\Vert
(e_{i}^{\ast })\right\Vert _{l_{s}^{n_{1}w}}.$

\noindent Then, $T$ is Lipschitz $\left( p,r,s\right) $-summing and we have%
\begin{equation*}
\pi _{p,r,s}^{SL}\left( T\right) \leq \pi _{p,r,s}\left( \widehat{T}\right)
.\quad \blacksquare 
\vspace{0.5cm}%
\end{equation*}

The following theorem gives an integral characterization of the strictly
Lipschitz $\left( p,r,s\right) $-summing operators. The proof is an
adaptation of the one in the last section, for this we will omit it.%
\vspace{0.5cm}%

\textbf{Theorem 4.5}. \textit{Suppose that }$\frac{1}{p}=\frac{1}{r}+\frac{1%
}{s}.$\textit{\ Let }$X$\textit{\ be a pointed metric space and }$E$\textit{%
\ be a Banach space. Let }$T:X\rightarrow E$\textit{\ be a Lipschitz
operator. The following properties are equivalent.}

\noindent \textit{1) }$T$\textit{\ is strictly Lipschitz }$\left(
p,r,s\right) $\textit{-summing.}

\noindent \textit{2) There exist a constant }$C>0$ \textit{and} \textit{a} 
\textit{Radon} \textit{probability }$\mu $\textit{\ on }$B_{X^{\#}}$ and $%
\eta \in B_{E^{\ast \ast }}$ \textit{such that for all }$x^{j},y^{j}\in X,$ $%
\lambda ^{j}\in \mathbb{K}\left( 1\leq j\leq n\right) $and $e^{\ast }\in
E^{\ast },$ \textit{we have}%
\begin{equation*}
\left\vert \left\langle \sum_{j=1}^{n}\lambda ^{j}\left( T\left(
x^{j}\right) -T\left( y^{j}\right) \right) ,e^{\ast }\right\rangle
\right\vert \leq C(\dint\limits_{B_{X^{\#}}}\left\vert \sum_{j=1}^{n}\lambda
^{j}\left( f\left( x^{j}\right) -f\left( y^{j}\right) \right) \right\vert
^{r}d\mu \left( f\right) )^{\frac{1}{r}}\left\Vert e^{\ast }\right\Vert
_{L_{s}\left( \eta \right) }.%
\vspace{0.5cm}%
\end{equation*}

\textbf{Corollary 4.6}. \textit{Let} $p\in \left[ 0,\infty \right] $ and $%
r,s\geq 1$. \textit{Suppose that }$\frac{1}{p}=\frac{1}{r}+\frac{1}{s}.$%
\textit{\ Let }$X$\textit{\ be a pointed metric space and }$E$\textit{\ be a
Banach space. Let }$T:X\rightarrow E$\textit{\ be a Lipschitz operator. The
following properties are equivalent.}

\noindent \textit{1) }$T$\textit{\ is strictly Lipschitz }$\left(
p,r,s\right) $\textit{-summing.}

\noindent \textit{2) There exist a Banach space }$G,$\textit{\ a strictly
Lipschitz }$s^{\ast }$\textit{-summing operator }$R:X\rightarrow G$\textit{\
and a Cohen strongly }$r$\textit{-summing linear operator }$S:G\rightarrow E$
\textit{such that}%
\begin{equation*}
T=S\circ R.%
\vspace{0.5cm}%
\end{equation*}

Let $1\leq p,r,s<\infty $ such that $\frac{1}{p}+\frac{1}{r}+\frac{1}{s}=1.$
Let $E,F$ be two Banach spaces, Lapreste \cite{13} has defined the norm $\mu
_{p,r,s}$ as follows%
\begin{equation*}
\mu _{p,r,s}\left( u\right) =\inf \left\{ \left\Vert \left( \lambda
_{i}\right) \right\Vert _{l_{p}^{n}}\left\Vert \left( e_{i}\right)
\right\Vert _{l_{r}^{n,w}\left( E\right) }\left\Vert \left( y_{i}\right)
\right\Vert _{l_{s}^{n,w}\left( F\right) }\right\} ,
\end{equation*}%
where the infimum is taken over all representations of $u=$ $%
\sum_{i=1}^{n}\lambda _{i}e_{i}\otimes y_{i}\in E\otimes F.$ For the
Lipschitz case, Ch\'{a}vez-Dom\'{\i}nguez in \cite{3}, has defined the
corresponding Lipschitz norm $\mu _{p,r,s}$\ (with the same notation of
Lapreste) on the space of molecules $\mathcal{F}\left( X;E\right) $ and
shown that the class $\Pi _{p^{\ast },r,s}^{L}\left( X,E^{\ast }\right) $
coincides with the dual of $\left( \mathcal{F}\left( X;E^{\ast }\right) ,\mu
_{p,r,s}\right) .$ Note that the space $\mathcal{F}\left( X;E\right) $ plays
the same role of Lipschitz tensor product $X\boxtimes E.$ Let $%
u=\sum_{i=1}^{n}\delta _{\left( x_{i},y_{i}\right) }\boxtimes e_{i}\in
X\boxtimes E^{\ast }$, the norm $\mu _{p,r,s}$ is defined as follows%
\begin{equation*}
\mu _{p,r,s}\left( u\right) =\inf \left\{ \left\Vert \left( \lambda
_{i}\right) \right\Vert _{l_{p}^{n}}\left\Vert \left( \lambda
_{i}^{-1}\kappa _{i}^{-1}\delta _{\left( x_{i},y_{i}\right) }\right)
\right\Vert _{l_{r}^{n,w}\left( \mathcal{F}\left( X\right) \right)
}\left\Vert \left( \kappa _{i}e_{i}^{\ast }\right) \right\Vert
_{l_{s}^{n,w}\left( E^{\ast }\right) }\right\} ,
\end{equation*}%
where the infimum is taken over all representations of $u$ in $X\boxtimes
E^{\ast }$ and $\lambda _{i},\kappa _{i}>0$. Let $X$ be a pointed metric
space and $E$ be a Banach space. Let $u=\sum_{k=1}^{l}\delta _{\left(
x_{k},y_{k}\right) }\boxtimes s_{k}^{\ast }\in X\boxtimes E^{\ast }$ and 
\begin{equation*}
A_{u}=\left\{ m=\sum_{i=1}^{n}\lambda _{i}m_{i}\otimes e_{i}^{\ast }\in 
\mathcal{F}\left( X\right) \otimes E^{\ast }:m=\sum_{k=1}^{l}\delta _{\left(
x_{k},y_{k}\right) }\otimes s_{k}^{\ast }\right\} .
\end{equation*}%
We consider%
\begin{equation*}
\mu _{p,r,s}^{SL}\left( u\right) =\inf_{m\in A_{u}}\left\{ \left\Vert \left(
\lambda _{i}\right) \right\Vert _{l_{p}^{n_{1}}}\left\Vert \left(
m_{i}\right) \right\Vert _{l_{r}^{n_{1},w}\left( \mathcal{F}\left( X\right)
\right) }\left\Vert \left( e_{i}^{\ast }\right) \right\Vert
_{l_{s}^{n_{1},w}\left( E^{\ast }\right) }\right\} .
\end{equation*}

In fact, in the Lipschitz definition of the norm $\mu _{p,r,s}$ we have just
used elements in $X\boxtimes E^{\ast }$ of the form $\sum_{i=1}^{n}\delta
_{\left( x_{i},y_{i}\right) }\boxtimes e_{i}^{\ast }$ which equals to $u$
and $\lambda _{i},\kappa _{i}>0,$ but in the definition of $\mu
_{p,r,s}^{SL},$ we have to consider all elements of the set $A_{u}.$
Therefore, the infimum in $\mu _{p,r,s}^{SL}$ will in general be smaller.$%
\vspace{0.5cm}%
$

\textbf{Proposition 4.7}. \textit{Let }$X$\textit{\ be a pointed metric
space and }$E$\textit{\ be a Banach space. For every }$u\in X\boxtimes E$%
\textit{\ we have}%
\begin{equation*}
\mu _{p,r,s}^{SL}\left( u\right) =\mu _{p,r,s}^{L}\left( u\right) ,
\end{equation*}%
\textit{where }$\mu _{p,r,s}^{L}$\textit{\ is the Lipschitz cross-norms
corresponding to the Lapreste tensor norm }$\mu _{p,r,s}.$%
\vspace{0.5cm}%

\textit{Proof}. Let $u=\sum_{k=1}^{l}\delta _{\left( x_{k},y_{k}\right)
}\boxtimes s_{k}^{\ast }\in X\boxtimes E^{\ast }.$ Let $m=\dsum%
\limits_{i}^{n}\lambda _{i}m_{i}\otimes e_{i}^{\ast }\in A_{u},$ we have

\QTP{Body Math}
$\mu _{p,r,s}^{SL}\left( u\right) =\inf_{m\in A_{u}}\left\{ \left\Vert
\left( \lambda _{i}\right) \right\Vert _{l_{p}^{n_{1}}}\left\Vert \left(
m_{i}\right) \right\Vert _{l_{r}^{n_{1},w}\left( \mathcal{F}\left( X\right)
\right) }\left\Vert \left( e_{i}^{\ast }\right) \right\Vert
_{l_{s}^{n_{1},w}\left( Y^{\ast }\right) }\right\} $

\QTP{Body Math}
$=\mu _{p,r,s}(\sum_{k=1}^{l}\delta _{\left( x_{k},y_{k}\right) }\otimes
s_{k}^{\ast })($by (2.3)$)$

\QTP{Body Math}
$=\mu _{p,r,s}^{L}\left( u\right) .\quad \blacksquare 
\vspace{0.5cm}%
$

\textbf{Corollary 4.8}. \textit{Let} $1\leq p\leq \infty .$ \textit{The
following properties are equvalent.}

\noindent \textit{1) The Lipschitz operator }$T:X\rightarrow E$\textit{\
strictly Lipschitz }$\left( p,r,s\right) $\textit{-summing.}

\noindent \textit{2) There exists a positive constant }$C$\textit{\ such
that for every }$x_{k},y_{k}\in X$\textit{\ and }$s_{k}^{\ast }\in E^{\ast }$%
\textit{\ }$\left( 1\leq k\leq l\right) $\textit{\ we have }%
\begin{equation*}
\left\vert \sum_{k=1}^{l}\left\langle T\left( x_{k}\right) -T\left(
y_{k}\right) ,s_{k}^{\ast }\right\rangle \right\vert \leq C\mu
_{p,r,s}^{SL}(u),
\end{equation*}%
\textit{where} $u=\sum_{k=1}^{l}\delta _{\left( x_{k},y_{k}\right)
}\boxtimes s_{k}^{\ast }.%
\vspace{0.5cm}%
$

Then, by \cite[Corollary 3.2]{16} we get the following coincidence result.$%
\vspace{0.5cm}%
$

\textbf{Corollary 4.9}. \textit{Let }$X$\textit{\ be a pointed metric space
and }$E$\textit{\ be a Banach space. Then}%
\begin{equation*}
\Pi _{p^{\ast },r,s}^{SL}\left( X,E^{\ast }\right) =\left( X\widehat{%
\boxtimes }_{\mu _{p,r,s}^{SL}}E^{\ast }\right) ^{\ast }=\left( \mathcal{F}%
\left( X\right) \widehat{\otimes }_{\mu _{p,r,s}}E^{\ast }\right) ^{\ast }.
\end{equation*}

\end{document}